\documentclass[11pt]{article}

\usepackage{blindtext}
\usepackage[preprint]{jmlr2e}
\usepackage[final]{changes}

\usepackage{url}
\usepackage{graphicx} 
\usepackage{multirow}
\usepackage{booktabs} 
\usepackage{pgfplotstable}

\usepackage{bm}
\usepackage{array} 
\usepackage{paralist} 
\usepackage{verbatim} 
\usepackage{appendix}
\usepackage{epstopdf}
\usepackage{bbm}
\usepackage{algorithm}
\usepackage{algorithmicx}
\usepackage{algpseudocode}
\usepackage{amsmath,nccmath,dsfont}
\usepackage{harpoon}
\usepackage{amssymb}
\usepackage{mathrsfs}
\usepackage{graphics}
\usepackage{color}
\usepackage{geometry}
\usepackage{epsfig}
\usepackage{diagbox}
\usepackage{mathtools}
\usepackage{bbm}
\usepackage{lineno}
\usepackage{fancyhdr}
\usepackage{graphicx}

\def\d{\text{d}}
\def\E{\mathbb{E}}
\usepackage{wrapfig}
\usepackage{dsfont}
\usepackage{enumerate}
\usepackage{subfigure}
\usepackage{xcolor}
\usepackage{setspace}

\usepackage{jmlr2e}

\usepackage{lastpage}
\jmlrheading{xx}{2024}{1-\pageref{LastPage}}{01/20; Revised xx/xx}{xx/xx}{21-0000}{Mingtao Xia and Xiangting Li and Qijing Shen and Tom Chou}

\ShortHeadings{Squared Wasserstein Distance for SDE
  Reconstruction}{Xia and Li and Shen and Chou} \firstpageno{1}

\begin{document}

\title{Squared Wasserstein-2 Distance for Efficient Reconstruction of Stochastic Differential Equations}
\author{\name Mingtao Xia\thanks{ equal contribution} \email xiamingtao@nyu.edu \\
\addr Courant Institute of Mathematical Sciences\\
New York University\\
New York, NY 10012, USA 
\AND
\name Xiangting Li$^*$ \email xiangting.li@ucla.edu \\
\addr Department of Computational Medicine \\
University of California Los Angeles \\
Los Angeles, CA 90095, USA
\AND
\name Qijing Shen \email qijing.shen@ndm.oxford.edu\\
\addr Nuffield Department of Medicine \\
University of Oxford \\
Oxford OX2 6HW, UK
\AND
\name Tom Chou \email tomchou@ucla.edu\\
\addr Department of Mathematics, \\
University of California Los Angeles\\
Los Angeles, CA 90095, USA}

\editor{TBD}

\maketitle

\begin{abstract}
We provide an analysis of the squared Wasserstein-2 ($W_2$) distance
between two probability distributions associated with two stochastic
differential equations (SDEs). Based on this analysis, we propose the
use of a squared $W_2$ distance-based loss functions in the
\textit{reconstruction} of SDEs from noisy data. To demonstrate the
practicality of our Wasserstein distance-based loss functions, we
performed numerical experiments that demonstrate the efficiency of
our method in reconstructing SDEs that arise across a number of
applications.
\end{abstract}

\begin{keywords}
  Wasserstein distance, stochastic differential equation, inverse
  problem, uncertainty quantification, optimal transport
\end{keywords}

\section{Introduction}
Stochastic processes are mathematical models of random phenomena that
evolve over time or space \citep{cinlar2011probability}. Among
stochastic processes, stochastic differential equations (SDE) of the
form
\begin{equation}
\d X(t) = f(X(t), t)\d t+\sigma(X(t), t)\d B(t), \,\,\, t\in[0, T]
\label{SDE_representation}
\end{equation}
are widely used across different fields to model complex systems with
continuous variables and noise. Here, $f$ and $\sigma$ denote
deterministic and stochastic components of the SDE, while $B(t)$
represents Brownian motion. In applications such as computational
fluid dynamics, cell biology, and genetics, the underlying dynamics
are often unknown, partially observed, and subjected to
noise. Consequently, it is vital to develop methods capable of
reconstructing the governing SDEs from limited data
\citep{sullivan2015introduction, soize2017uncertainty,
  mathelin2005stochastic, bressloff2014stochastic,
  lin2018efficient}. Traditional methods, such as the Kalman filter
\citep{welch1995introduction, welch2020kalman} and Gaussian process
regression \citep{liu2020gaussian,mackay1998introduction} often assume
specific forms of noise. These methods may not be suitable for complex
or nonlinear systems where noise affects the dynamics in a more
complex manner.

Recent advancements leverage machine learning, specifically neural
ordinary differential equations (NODEs) \citep{chen2018neural}, to
offer a more flexible approach to reconstructing SDEs in the form of
neural SDEs (nSDEs)
\citep{tzen2019neural,tong2022learning,jia2019neural}. Despite the
promise, challenges remain, particularly in selecting optimal loss
functions \citep{jia2019neural}. The Wasserstein distance, a family of
metrics that measures discrepancies between probability measures over
a metric space, has emerged as a potential solution due to its robust
properties \citep{villani2009optimal,
  oh2019kernel,zheng2020nonparametric}. In this paper, we introduce
bounds on the second-order Wasserstein $W_2$ distance between two
probability distributions over the continuous function space generated
by solutions to two SDEs. Our results motivate the use of this
distance for SDE reconstruction. We test our approach on different
examples to showcase its effectiveness.

Traditional methods for reconstructing SDEs from data usually make
assumptions on the specific forms of the underlying SDE and fit
unknown parameters. For example, \citep{de2016reduction} uses some
polynomials to model $f, \sigma$, while \citep{pereira2010learning}
assumes linear $f$ and $\sigma$ in Eq.~\eqref{SDE_representation}.

Previous attempts at using neural SDEs (nSDEs) have explored different
loss functions for reconstruction. For example, \cite{tzen2019neural}
model the SDE as a continuum limit of latent deep Gaussian models and
use a variational likelihood bound for
training. \cite{kidger2021neural} adopt Wasserstein generative
adversarial networks (WGANs) that were proposed in
\cite{arjovsky2017wasserstein} for reconstructing
SDEs. \cite{briol2019statistical} uses a maximum mean discrepancy
(MMD) loss and a generative model for training
SDEs. \cite{song2020score} assumes that $\sigma$ in
Eq.~\eqref{SDE_representation} depends only on time and uses a
score-based generative model for SDE reconstruction.

The Wasserstein distance, denoted as $W$, has gained wide use in
statistics and machine learning. Key papers have delved into its
analysis \citep{ruschendorf1985wasserstein} and its utilization in
reconstructing discrete-time stochastic processes
\citep{bartl2021wasserstein}. In the context of SDEs,
\cite{bion2019wasserstein} introduced a restricted Wasserstein-type
distance, while \cite{sanz2021wasserstein} and \cite{wang2016p,
  sanz2021wasserstein} examined its application in ergodic SDEs, Levy
processes, and Langevin equations, respectively. Calculating the $W$
distance for multidimensional random variables is challenging; hence,
approximations such as the sliced \( W \) distance and regularized \(
W \) distance have emerged \citep{cuturi2019differentiable,
  kolouri2018sliced, kolouri2019generalized, rowland2019orthogonal,
  frogner2015learning}.

The aforementioned WGAN approach in \cite{kidger2021neural} uses the
first-order Wasserstein distance to indirectly reconstruct SDEs via
the Kantorovich-Rubinstein duality \citep{arjovsky2017wasserstein}. To
the best of our knowledge, there has been no published work that
directly applies and analyzes the $W$ distance to the reconstruction
of SDEs.

\section{Definitions and Outline}
We propose a squared $W_2$-distance-based SDE reconstruction method
and analyze it under the following setting. Let $\mu$ denote the
probability distribution over the continuous function space $C([0, T];
\mathbb{R})$ generated by the solution $X(t)$ to
Eq.~\eqref{SDE_representation}. In the following approximation to
Eq.~\eqref{SDE_representation},
\begin{equation}
  \text{d} \hat{X}(t) = \hat{f}(X(t), t)\d t
  + \hat{\sigma}(X(t), t)\text{d} \hat{B}(t), \,\,\, t\in[0, T],
\label{approximate_sde}
\end{equation}
$\hat{B}(t)$ is another independent standard Brownian motion and the
probability distribution over the continuous function space $C([0, T];
\mathbb{R})$ generated by the solution $\hat{X}(t)$ to
Eq.~\eqref{approximate_sde} will be denoted $\hat{\mu}$.

We shall follow the definition of the squared $W_2$-distance in
\cite{clement2008elementary} for two probability measures $\mu,
\hat{\mu}$ associated with two continuous stochastic processes $X(t),
\hat{X}(t),\,\, t\in[0, T]$.

\begin{definition} 
\rm
\label{def:W2}
For two $d$-dimensional continuous stochastic processes in the
separable space $\big(C([0, T]; \mathbb{R}^d), \|\cdot\|\big)$
\begin{equation}
\bm{X}(t)=(X^1(t),...,X^d(t)), \,\,\hat{\bm{X}}(t)
= (\hat{{X}}^1(t),...,\hat{X}^d(t)),\, t\in[0, T],
\label{sde_dimension}
\end{equation}
with two associated probability distributions $\mu, \hat{\mu}$, the
squared $W_2(\mu, \hat{\mu})$ distance between $\mu, \hat{\mu}$ is
defined as
\begin{equation}
W_2^2(\mu, \hat{\mu}) = \inf_{\pi(\mu, \hat \mu)}
\E_{(\bm{X}, \hat{\bm{X}})\sim \pi(\mu, \hat \mu)}\big[\|{\bm{X}} - \hat{{\bm{X}}}\|^2\big].
\label{pidef}
\end{equation}
The distance $\|\bm{X}\|\coloneqq \big(\int_0^T \sum_{i=1}^d
|X_i(t)|^2\d t\big)^{\frac{1}{2}}$ and $\pi(\mu, \hat \mu)$ iterates
over all \textit{coupled} distributions of $\bm{X}(t),
\hat{\bm{X}}(t)$, defined by the condition

\begin{equation}
\begin{cases}
  {\bm{P}}_{\pi(\mu, \hat \mu)}\left(A \times C([0, T]; \mathbb{R}^d)\right)
  ={\bm{P}}_{\mu}(A),\\
  {\bm{P}}_{\pi(\mu, \hat \mu)}\left(C([0, T]; \mathbb{R}^d)\times A\right)
  = {\bm{P}}_{\hat \mu}(A), 
\end{cases}\forall A\in \mathcal{B}\Big(C([0, T]; \mathbb{R}^d)\Big),
\end{equation}
where $\mathcal{B}\Big(C([0, T]; \mathbb{R}^d)\Big)$ denotes the Borel
$\sigma$-algebra associated with the space of $d$-dimensional
continuous functions $C([0, T]; \mathbb{R}^d)$.
\end{definition}

\noindent Our main contributions can be summarized as follows
\begin{itemize}
\item[$1.$] Using Definition \ref{def:W2}, we first derive in
  Section~\ref{section2} an upper bound for the squared Wasserstein
  distance $W_2^2(\mu, \hat{\mu})$ between the probability measures
  associated with solutions to two 1D SDEs in terms of the errors in the
  reconstructed drift and diffusion functions $f-\hat{f}$ and
  $\sigma-\hat{\sigma}$ in Eqs.~\eqref{SDE_representation} and
  \eqref{approximate_sde}.  To be specific, we establish a $W_2$
  distance upper bound which depends explicitly on the difference in
  the drift and diffusion functions $f-\hat{f}$ and
  $\sigma-\hat{\sigma}$ associated with using
  Eq.~\eqref{approximate_sde} to approximate
  Eq.~\eqref{SDE_representation}.
\item[$2.$] In Section \ref{section25}, we shall prove that the $W_2$
  distance between the two SDEs, $W_2(\mu, \hat{\mu})$, can be
  accurately approximated by estimating the $W_2$ distance between
  their finite-dimensional projections.  We also develop a
  time-decoupled version of $W_2^2(\mu, \hat{\mu})$ defined by
\begin{equation}
\tilde{W}_2^2(\mu,\hat{\mu})\coloneqq \int_0^TW_2^2(\mu(s), \hat{\mu}(s))\text{d}s
\label{time_decoupled_W2}
\end{equation}
which allows us to define a squared $W_2$-distance-based loss function
for reconstructing SDEs.  Here, $\mu(s), \hat{\mu}(s)$ are the
distributions on $\mathbb{R}^d$ generated by projection of the
stochastic processes $\bm{X}, \hat{\bm{X}}$ at time $s$, respectively.
Specifically, if $X(t_i)$ follows a one-dimensional SDE
Eq.~\eqref{SDE_representation}, then for uniformly spaced time points
$t_i=\frac{iT}{N},\,\, i=0,...,N$, our proposed time-decoupled loss
function is simply
\vspace{-0.05in}
\begin{equation}
\Delta{t}\sum_{i=1}^{N-1}\int_0^1 (F_{i}^{-1}(s) - \hat{F}_{i}^{-1}(s))^2\text{d}s ,
\label{new_loss}
\end{equation}
where $\Delta{t}$ is the timestep and $F_{i}$ and $\hat{F}_{i}$ are
the empirical cumulative distribution functions for $X(t_i)$ and
$\hat{X}(t_i)$, respectively. This time-decoupled squared
$W_2$-distance loss function will be explicitly expressed in
Eq.~\eqref{approximation}.
\item[$3.$] Finally, we carry out numerical experiments to show that our
  squared $W_2$-distance-based SDE reconstruction method performs
  better than recently developed machine-learning-based methods across
  many SDE reconstruction problems. In Section~\ref{summary}, we
  summarize our proposed squared $W_2$ distance method for SDE
  reconstruction and suggest some promising future directions.
  Additional numerical experiments and sensitivity analysis are
  detailed in the Appendix.
\end{itemize}

\section{Squared $W_2$ distance for reconstructing SDEs}
\label{section2}

In this section, we prove the bounds for the squared $W_2$ distance of
two probability measures associated with two SDEs. Specifically, we
demonstrate that minimizing the squared $W_2$ distance is necessary
for the reconstruction of $f, \sigma$ in
Eq.~\eqref{SDE_representation}.

We shall first prove an upper bound for the $W_2$ distance between the
probability measures $\mu$ and $\hat\mu$ associated with $X(t),
\hat{X}(t)$, solutions to Eq.~\eqref{SDE_representation} and
Eq.~\eqref{approximate_sde}, respectively.
\setcounter{theorem}{0}
\begin{theorem}
\label{theorem1}
\rm If $\{X(t)\}_{t=0}^T, \{\hat{X}(t)\}_{t=0}^T$ have the same
initial condition distribution and they are solutions to
Eq.~\eqref{SDE_representation} and Eq.~\eqref{approximate_sde}
in the univariate case ($d=1$ in Eq.~\eqref{sde_dimension}),
respectively, and the following conditions hold:

\begin{itemize}
\item $f, \hat{f}, \sigma, \hat{\sigma}$ are continuously
  differentiable; $\partial_x\sigma$ and $\partial_x{\hat{\sigma}}$
  are uniformly bounded
\item there exists two functions $\eta_1(x_1, x_2), \eta_{2}(x_1,
  x_2)$ such that their values are in $(x_1, x_2)$ and
\begin{equation}
\begin{aligned}
&f(X_1, t) - f(X_2, t) = \partial_{x}f(\eta_{1}(X_1, X_2), t) (X_1-X_2)\,\,\,\\
&\sigma(X_1, t) -
\sigma(X_2, t) = \partial_{x} \sigma(\eta_{2}(X_1, X_2),
t)(X_1-X_2) 
\end{aligned}
\label{intermediate_eq}
\end{equation}
\end{itemize}
then,

\begin{small}
\begin{equation}
\begin{aligned}
W_2^2(\mu, \hat{\mu}) \leq & 3\int_0^T
\E\Big[\medint\int_0^t H^2(s, t)\d s \Big]\d t \times
\E\Big[\medint\int_0^T (f- \hat{f})^2(\tilde{X}(t), t)\d t)\Big]\\
\: & +3\int_0^T \!\E\Big[\medint\int_0^t \!H^2(s, t)\d s\Big]\d t \times
\E\Big[\medint\int_0^T
\!\big(\partial_x\sigma(\eta_2 (X(t), \tilde{X}(t)), t\big)^2
(\sigma- \hat{\sigma})^2(\tilde{X}(t), t)\d t\Big] \\ 
\: & +3\int_0^T \!\E\Big[\medint\int_0^t \!H^4(s, t)\d s\Big]^{1/2}\d t
\times \E\Big[\medint\int_0^T
  \!(\sigma- \hat{\sigma})^4(\tilde{X}(t), t)\d t)\Big]^{1/2},
\end{aligned}
\label{W2_estimate}
\end{equation}
\end{small}
where $\tilde{X}(t)$ satisfies

\begin{equation}
\text{d}\tilde{X}(t) = \hat{f}(\tilde{X}(t), t)\d t
+ \hat{\sigma}(\tilde{X}(t), t)\text{d} B(t),\,\,\, \tilde{X}(0) = X(0),
\label{tilde_x}
\end{equation}
and

\begin{equation}
H(s,t)\coloneqq \exp
\left[{\int_s^t h(X(r), \tilde{X}(r),r)\d r + \int_s^t
\partial_x\sigma\big(\eta_2 (X(r),\tilde{X}(r), r\big) \d
B(r)}\right],
\end{equation}
with $h$ defined as

\begin{equation}
\begin{aligned}
  h(X(r), \hat{X}(r), r) \coloneqq \partial_x
  f\big(\eta_1(X(r),\tilde{X}(r)), r\big)
-\Big(\partial_x\sigma\big(\eta_2(X(r),\tilde{X}(r)),r\big)\Big)^2.
\end{aligned}
\label{h_def}
\end{equation}
\end{theorem}

The proof to Theorem~\ref{theorem1} and its generalizations to higher
dimensional stochastic dynamics under some specific assumptions are given in
Appendix~\ref{proof_theorem1}.  Theorem~\ref{theorem1} indicates that
as long as $\E\big[\int_0^t H^{4}(s, t) \d s\big]$ is uniformly
bounded for all $t\in[0, T]$, the upper bound for $W_2(\mu,
\hat{\mu})\rightarrow 0$ when $\hat{f}-f\rightarrow 0$ and
$\hat{\sigma}-\sigma\rightarrow 0$ uniformly in $\mathbb{R}\times[0,
  T]$.
Specifically, if $f=\hat{f},\sigma=\hat{\sigma}$, then the RHS
Eq.~\eqref{W2_estimate} is 0. This means that minimizing $W_2^2(\mu,
\hat{\mu})$ is necessary for generating small errors $\hat{f}-f,
\hat{\sigma}-\sigma$ and for accurately approximating both $f$ and
$\sigma$. Thus, one can consider using the squared $W_2$ distance as
an effective loss function to minimize during reconstruction of SDEs.
MSE-based loss functions (defined in Appendix \ref{def_loss}) suppress
noise while the Kullback-Liebler (KL) divergence may not be finite,
precluding resolution of \(X(t)\) and \(\hat{X}(t)\) even if
\(\hat{f}\) approximates \(f\) and \(\hat{\sigma}\) approximates
\(\sigma\).  Detailed discussions on the limitations of MSE and KL
divergence in SDE reconstruction can be found in
Appendix~\ref{proof_theorem2}.

\section{Finite-dimensional and time-decoupled squared $W_2$ loss functions}
\label{section25}

From Theorem~\ref{theorem1} in Section~\ref{section2}, in order to
have small errors in the drift and diffusion terms $f-\hat{f}$ and
$\sigma-\hat{\sigma}$, a small $W_2(\mu, \hat{\mu})$ is
necessary. However, $W_2(\mu, \hat{\mu})$ cannot be directly used as a
loss function to minimize since we cannot directly evaluate the
integration in time in Eq.~\eqref{pidef}.  In this section, we shall
provide a way to estimate the $W_2(\mu, \hat{\mu})$ distance by using
finite dimensional projections, leading to squared
$W_2$-distance-based loss functions for minimization.

Consider the two general $d$-dimensional SDEs defined in
Eq.~\eqref{multi_dimensional}.  Usually, we only have finite
observations of trajectories for $\{\bm{X}(t)\}$ and
$\{\hat{\bm{X}}(t)\}$ at discrete time points. Thus, we provide an
estimate of the $W_2$ between of the probability measures $\mu,
\hat{\mu}$ associated with $\bm{X}(t)$ and $\bm{\hat{X}}(t)$, $t\in[0,
  T]$ using their finite-dimensional projections. We assume that
$\bm{X}(t), \bm{\hat{X}}(t)$ solve the two SDEs described by
Eq.~\eqref{multi_dimensional}. We let $0=t_0<t_1<...<t_N=T,
t_i=i\Delta{t}, \Delta{t}\coloneqq\frac{T}{N}$ be a uniform mesh in
time and we define the following projection operator $\bm{I}_N$
\begin{equation}
\bm{X}_N(t) \coloneqq \bm{I}_N \bm{X}(t) =\left\{
\begin{aligned}
&\bm{X}(t_i), t\in[t_i, t_{i+1}), i<N-1,\\
&\bm{X}(t_i), t\in[t_i, t_{i+1}], i=N-1.
\end{aligned}
\right.
\label{X_N_def}
\end{equation}
As in the previous case, we require $\boldsymbol{X}(t)$ and
$\hat{\boldsymbol{X}}(t)$ to be continuous. Note that the projected
process is no longer continuous. Thus, we define a new space
$\tilde\Omega_{N}$ containing all continuous and piecewise constant
functions; naturally, $\mu, \hat{\mu}$ are allowed to be defined on
$\tilde \Omega_N$.  Distributions over $\tilde \Omega_N$ generated by
$\boldsymbol{X}_N(t), \hat{\boldsymbol{X}}_N(t)$ in
Eq. \eqref{X_N_def} is denoted by $\mu_N$ and $\hat{\mu}_N$,
respectively. We will prove the following theorem for
  estimating $W_2(\mu, \hat{\mu})$ by $W_2\left(\mu_N,
  \hat{\mu}_N\right)$.

\begin{theorem}
\label{theorem3}
\rm Suppose $\{\bm{X}(t)\}_{t=0}^T$ and $\{\hat{\bm{X}}(t)\}_{t=0}^T$
are both continuous-time continuous-space stochastic processes in
$\mathbb{R}^{d}$ and $\mu, \hat{\mu}$ are their associated probability
measures, then $W_2(\mu, \hat{\mu})$ can be bounded by their
finite-dimensional projections

\begin{equation}
W_2(\mu_N, \hat{\mu}_N) - W_2(\mu, \mu_N)
- W_2(\hat{\mu}, \hat{\mu}_N)
\leq W_2(\mu, \hat{\mu})
\leq W_2(\mu_N, \hat{\mu}_N) + W_2(\mu, \mu_N) + W_2(\hat{\mu}, \hat{\mu}_N)
\label{triangular}
\end{equation}
where $\mu_N, \hat{\mu}_N$ are the probability distributions
associated with $\bm{X}_N$ and $\hat{\bm{X}}_N$ defined in
Eq.~\eqref{X_N_def}. Specifically, if $\bm{X}(t)$ and
$\hat{\bm{X}(t)}$ solve Eq.~\eqref{multi_dimensional}, and if
\begin{equation}
\begin{aligned}
&F\coloneqq\E\Big[\medint\int_0^T \sum_{i=1}^d f_i^2(\bm{X}(t),t)\d t\Big]<\infty,
\,\,\Sigma\coloneqq\E\Big[\medint\int_0^T \sum_{\ell=1}^d
\sum_{j=1}^s\sigma_{i, j}^2(\bm{X}(t),t)\d t\Big]<\infty,\\
&\hat{F}\coloneqq\E\Big[\medint\int_0^T \sum_{i=1}^d
\hat{f}_i^2(\hat{\bm{X}}(t),t)\d t\Big]<\infty,\,\,
\hat{\Sigma}\coloneqq\E\Big[\medint\int_0^T
\sum_{\ell=1}^d\sum_{j=1}^s\hat{\sigma}_{i, j}^2(\hat{\bm{X}}(t),t)\d t\Big]<\infty,
\end{aligned}
\label{F_Sigma}
\end{equation}
then we obtain the following bound
\begin{equation}
\begin{aligned}
W_2(\mu_N, \hat{\mu}_N)- & \sqrt{(s+1)\Delta t}
\left(\sqrt{F\Delta t+\Sigma}+\sqrt{\hat{F}\Delta t
+\hat{\Sigma}}\right) \leq W_2(\mu,\hat{\mu}) \\
\: &  \leq W_2(\mu_N, \hat{\mu}_N)
+\sqrt{(s+1)\Delta t}\left(\sqrt{F\Delta t+\Sigma}
+\sqrt{\hat{F}\Delta t+\hat{\Sigma}}\right).
\end{aligned}
\label{dtbound}
\end{equation}
\end{theorem}


The proof to Theorem~\ref{theorem3} relies on the triangular inquality
of the Wasserstein distance and the It\^o isometry; it is provided in
Appendix~\ref{proof_theorem3}. Theorem~\ref{theorem3} gives bounds for
approximating the $W_2$ distance between $\bm{X}(t), \hat{\bm{X}}(t)$
w.r.t. to their finite dimensional projections $\bm{X}_N(t),
\hat{\bm{X}}_N(t)$. Specifically, if $\bm{X}(t), \hat{\bm{X}}(t)$ are
solutions to Eq.~\eqref{SDE_representation} and
Eq.~\eqref{approximate_sde}, then as the timestep $\Delta t\rightarrow
0$, $W_2(\mu_N, \hat{\mu}_N)\rightarrow W_2(\mu, \hat{\mu})$.
Theorem~\ref{theorem3} indicates that we can use $W_2^2(\mu_N,
\hat{\mu}_N)$, which approximates $W_2^2(\mu, \hat{\mu})$ when
$\Delta{t}\rightarrow0$, as a loss function. Furthermore,

\begin{equation}
W_2^2(\mu_N, \hat{\mu}_N) = {\inf}_{\pi(\mu_N, \hat\mu_N)}
\sum_{i=1}^{N-1} \E_{(\bm X_N, \hat {\bm X}_N)\sim\pi(\mu_N,
  \hat{\mu}_N)}\Big[\big|\bm{X}(t_i) -
  \hat{\bm{X}}(t_i)\big|_2^2\Big]\Delta t.
\label{time_discretize}
\end{equation}
Here, $\pi(\mu_N, \hat\mu_N)$ iterates over coupled distributions of
$\bm{X}_N(t), \hat{\bm{X}}_N(t)$, whose marginal distributions
coincide with $\mu_N$ and $\hat{\mu}_N$. $|\cdot|_2$ denotes the
$\ell^2$ norm of a vector. Note that $\mu_N$ is fully characterized by
values of $\bm{X}(t)$ at the discrete time points $t_i$.

For a $d$-dimensional SDE, the trajectories at discrete time points
$\{\bm{X}(t_i)\}_{i=1}^{N-1}$ is $d\times (N-1)$ dimensional. In
\cite{fournier2015rate}, the error bound for $|W_2^2(\mu_N,
\hat{\mu}_N) - W_2^2(\mu^{\text{e}}_N, \hat{\mu}^{\text{e}}_N)|$,
where $\mu^{\text{e}}_N, \hat{\mu}^{\text{e}}_N$ are the finite-sample
empirical distributions of $\{\bm{X}(t_i)\}_{i=1}^{N-1}$ and
$\{\hat{\bm{X}}(t_i)\}_{i=1}^{N-1}$, will increase as the
dimensionality $d\times (N-1)$ becomes large.  Alternatively, we can
disregard the temporal correlations of values at different times and
relax the constraint on the coupling $\pi(\mu_N, \hat\mu_N)$ in to
minimize the Wasserstein distance between the marginal distribution of
$\{\bm{X}(t_i)\}$ and the marginal distribution of
$\{\hat{\bm{X}}(t_i)\}$, as was done in \cite{chewi2021fast}. To be
more specific, we minimize individual terms in the sum with respect to
the coupling $\pi_i$ of $\bm{X}(t_i)$ and $\hat{\bm{X}}(t_i)$ and
define a heuristic loss function

\begin{equation}
\sum_{i=1}^{N-1}
\inf_{\pi_i} \E_{\pi_i} \left[\big|\bm{X}(t_i) -
  \hat{\bm{X}}(t_i)\big|_2^2\right]
\Delta t =\sum_{i=1}^{N-1} W_2^2(\mu_N(t_i), \hat{\mu}_N(t_i))\Delta{t}
\label{approximation}
\end{equation}
where $\mu_N(t)$ and $\hat{\mu}_N(t)$ are the probability
distributions of $\bm{X}(t)$ and $\bm{X}(t)$, respectively.  Note that

\begin{equation}
\sum_{i=1}^{N-1}
\inf_{\pi_i} \E_{\pi_i} \big[\big|\bm{X}(t_i) -
\hat{\bm{X}}(t_i)\big|_2^2\big] \Delta t\leq W_2^2(\mu_N, \hat{\mu}_N) 
\end{equation}
because the marginal distributions of $\pi(\mu_N, \hat{\mu}_N)$
  coincide with $\mu_N$ and $\hat{\mu_N}$. Since the marginal distributions
of $\mu_N$ and $\hat{\mu}_N$ at $t_i$ are $\mu_N(t_i)$ and $\hat{\mu}_N(t_i)$,
respectively, we have

\begin{equation}
  \begin{aligned}
& \sum_{i=1}^{N-1}\inf_{\pi_i} \E_{\pi_i} \left[\big|\bm{X}(t_i) -
      \hat{\bm{X}}(t_i)\big|_2^2\right] \Delta t \\
& \qquad \qquad \leq {\inf}_{\pi(\mu_N, \hat\mu_N)} \sum_{i=1}^{N-1} 
\E_{(\bm X_N, \hat {\bm X}_N)\sim\pi(\mu_N, \hat{\mu}_N)}\Big[\big|\bm{X}(t_i)
  - \hat{\bm{X}}(t_i)\big|_2^2\Big]\Delta t.
\end{aligned}
\end{equation}

The dimensionality of $\bm{X}(t_i)$ and $\hat{\bm{X}}(t_i)$ is $d$,
which is much smaller than $(N-1)d$ for large $N$. We denote
$\mu^{\text{e}}_N(t_i)$ and $\hat{\mu}^{\text{e}}_N(t_i)$ to be the
finite-sample empirical distributions of $\bm{X}(t_i)$ and
$\hat{\bm{X}}(t_i)$, respectively.  Since the error of estimating the
$W_2$ distance using empirical distributions of a random variable
increases with the random variable's dimensionality
\cite{fournier2015rate}, the error
$\big|\sum_{i=1}^{N-1}W_2^2(\mu_N(t_i),
\hat{\mu}_N(t_i))-\sum_{i=1}^{N-1}W_2^2(\mu_N^{\text{e}}(t_i),
\hat{\mu}_N^{\text{e}}(t_i))\big|$ can be smaller than the error
$\big|W_2^2(\mu_N, \hat{\mu}_N) - W_2^2(\mu^{\text{e}}_N,
\hat{\mu}^{\text{e}}_N)\big|$. Compared to
Eq.~\eqref{time_discretize}, the time-decoupled squared $W_2$ distance
Eq.~\eqref{approximation} can be better approximated using
finite-sample empirical distributions.

Note that

\begin{equation}
  \sum_{i=1}^{N-1} W_2^2(\mu_N(t_i),
  \hat{\mu}_N(t_i))\Delta{t} \leq W_2^2(\mu_N, \hat{\mu}_N).
\end{equation}
Thus, from Theorem~\ref{theorem1} and Theorem~\ref{theorem3},
minimizing Eq.~\eqref{approximation} when $N\rightarrow \infty$ is
also necessary to achieve small $f-\hat{f}$ and $\sigma-\hat{\sigma}$
when the SDE is univariate.  Let $\boldsymbol{\mu}_i,
  \hat{\boldsymbol{\mu}}_i$ be the two probability distributions on
  the space of continuous functions associated with $\bm{X}(t),
  t\in[t_i, t_{i+1})$ and $\hat{\bm{X}}(t), t\in[t_i, t_{i+1})$, respectively. We
      can then show that Eq.~\eqref{approximation} is an approximation
      to the partially time-decoupled summation of squared $W_2$
      distances $\sum_{i=1}^{N-1}W_2^2(\boldsymbol{\mu}_i,
      \hat{\boldsymbol{\mu}}_i)$ as $N\rightarrow \infty$.
    Additionally, we can prove the following theorem that indicates
    Eq.~\eqref{approximation} approximates a time-decoupled squared
    Wasserstein distance in the $N\rightarrow\infty$ limit.

\begin{theorem}
\rm
\label{theorem4}

We assume the conditions in Theorem~\ref{theorem3} hold and for any
$0<t<t'<T$, as $t'-t\rightarrow 0$, the following conditions are
satisfied

\begin{equation}
\begin{aligned}
\E\Big[\medint\int_{t}^{t'} \sum_{i=1}^d f_i^2(\bm{X}(t),t)\text{d}t\Big],\,\,
  \E\Big[\medint\int_{t}^{t'} \sum_{i=1}^d \hat{f}_i^2(\hat{\bm{X}}(t),t)
    \text{d}t\Big] & \rightarrow 0,\\
  \E\Big[\medint\int_{t}^{t'} \sum_{i=1}^d
    \sum_{j=1}^s\sigma_{i, j}^2(\bm{X}(t),t)\d t\Big],\,\,
  \E\Big[\medint\int_{t}^{t'} \sum_{i=1}^d
\sum_{j=1}^s\hat{\sigma}_{i, j}^2(\hat{\bm{X}}(t),t)\d t\Big] & \rightarrow 0. 
\end{aligned}
\label{condition_dt}
\end{equation}
Then, 
\begin{equation}
  \lim\limits_{N\rightarrow\infty}\Big(\sum_{i=1}^{N-1}\inf_{\pi_i}
  \E_{\pi_i}\left[\big|\bm{X}(t_i) - \hat{\bm{X}}(t_i)\big|_2^2\right]\Delta t
  - \sum_{i=1}^{N-1}W_2(\boldsymbol{\mu}_i, \hat{\boldsymbol{\mu}}_i) \Big)= 0.
\end{equation}
Furthermore, the limit

\begin{equation}
  \lim\limits_{N\rightarrow \infty}\sum_{i=1}^{N-1}
  \inf_{\pi_i}\E_{\pi_i}\left[\big|\bm{X}(t_i) -
    \hat{\bm{X}}(t_i)\big|_2^2\right]\Delta t
  = \lim\limits_{N\rightarrow \infty}\sum_{i=1}^{N-1}
  W_2^2\Big(\mu(t_i), \hat{\mu}(t_i)\Big)\Delta t
\end{equation}
exists.

\end{theorem}

The proof of Theorem~\ref{theorem4} will use the result of
Theorem~\ref{theorem3} and is given in
Appendix~\ref{proof_theorem4}. Specifically, for each $N$,

\begin{equation}
  \sum_{i=1}^{N-1}\inf_{\pi_i}\E_{\pi_i}\left[\big|\bm{X}(t_i)
    - \hat{\bm{X}}(t_i)\big|_2^2\right]\Delta t \leq W_2^2(\mu_N, \hat{\mu}_N),
\end{equation}
so we conclude that

\begin{equation}
  \lim\limits_{N\rightarrow \infty} \sum_{i=1}^{N-1}\inf_{\pi_i}
  \E_{\pi_i}\left[\big|\bm{X}(t_i)
    - \hat{\bm{X}}(t_i)\big|_2^2\right]\Delta t
  \leq \lim\limits_{N\rightarrow\infty}W_2^2(\mu_N, \hat{\mu}_N)
  =W_2^2(\mu, \hat{\mu}).
\label{limit_relation}
\end{equation}
We denote

\begin{equation}
  \tilde{W}^2_2(\mu, \hat{\mu})\coloneqq
  \int_0^T W_2^2\big(\mu(t), \hat{\mu}(t)\big)\text{d}t
  =\lim\limits_{N\rightarrow \infty}\sum_{i=1}^{N-1}
  W_2^2\big(\mu(t_i^1), \hat{\mu}(t_i^1)\big)(t_i^1-t_{i-1}^1)
\end{equation}
as the \textit{time-decoupled squared Wasserstein distance}.  From
Eq.~\eqref{limit_relation}, we can deduce that

\begin{equation}
\tilde{W}_2^2(\mu, \hat{\mu}) \leq W_2^2(\mu, \hat{\mu}).
\end{equation}
Therefore, the upper bound of $W_2^2(\mu, \hat \mu)$ in
Theorem~\ref{theorem1} is also an upper bound of $\tilde W_2^2(\mu,
\hat \mu)$, \textit{i.e.}, to reconstruct a 1D SDE by
  minimizing $\tilde W_2^2(\mu, \hat \mu)$, it is necessary that
  $f-\hat{f}$ and $\sigma-\hat{\sigma}$ are small.  From
Theorem~\ref{theorem4}, minimizing the finite-time-point
time-decoupled loss function defined in Eq.~\eqref{approximation},
which approximates $\tilde{W}_2^2(\mu, \hat{\mu})$ when $\Delta t$ is
small, is needed for minimizing $f-\hat{f}$ and
$\sigma-\hat{\sigma}$.


Specifically, if $X(t), \hat{X}(t)$ are solutions to the univariate
SDEs Eq.~\eqref{SDE_representation} and Eq.~\eqref{approximate_sde},
then Eq.~\eqref{approximation} reduces to Eq.~\eqref{new_loss}, which
can be directly calculated.  In Example~\ref{example3}, Example
\ref{example4}, and Appendix~\ref{appendix_more_discussion}, we shall
compare use of the two different squared $W_2$ distance loss functions
Eqs.~\eqref{time_discretize} and ~\eqref{approximation}. From our
preliminary numerical results, using Eq.~\eqref{approximation} is more
efficient than using Eq.~\eqref{time_discretize} and gives more
accurate reconstructed SDEs.

\section{Numerical experiments}
\label{section3}
We carry out experiments to investigate the efficiency of our proposed
squared $W_2$ loss function (Eq.~\eqref{approximation}) by comparing
it to other methods and loss functions. Our approach is tested on the
reconstruction of several representative SDEs in
Examples~\ref{example1}--~\ref{example4}.

In all experiments, we use two neural networks to parameterize
$\hat{f}\coloneqq \hat{f}(X, t;\Theta_1), \hat{\sigma}\coloneqq
\hat{\sigma}(X, t;\Theta_2)$ in Eq.~\eqref{approximate_sde} for the
purpose of reconstructing $f, \sigma$ in
Eq.~\eqref{SDE_representation} by the estimates $\hat{f}\approx f,
\hat{\sigma}\approx{\sigma}$. $\Theta_1, \Theta_2$ are the parameter
sets in the two neural networks for parameterizing
$\hat{f}=\hat{f}_{\Theta_1}, \hat{\sigma}=\hat{\sigma}_{\Theta_2}$.
We use the \texttt{sdeint} function in the \texttt{torchsde}
Python package in \cite{li2020scalable} to numerically
integrate SDEs. Details of the training hyperparameter setting for all examples are
given in Appendix~\ref{training_details}.  Our code will be
  made publicly available on Github upon acceptance of this
  manuscript.

First, we compare our proposed squared $W_2$-distance-based loss
(Eq.~\eqref{approximation}) with several traditional statistical
methods for SDE reconstruction.

\begin{example}
\rm
\label{example1}
We reconstruct a nonlinear SDE of the form
\begin{equation}
  \text{d} X(t) = \big(\tfrac{1}{2}-\cos X(t)\big)\text{d} t
  + \sigma \text{d} B(t),\,\,\, t\in[0, 20],
\label{example1_ground_truth}
\end{equation}
which defines a Brownian process in a potential of the form $U(x) =
\frac{x}{2}-\sin x$. In the absence of noise, there are infinitely
many stable equilibrium points $x_k = \frac{5\pi}{3} + 2\pi k,
k\in\mathbb{Z}$. When noise $\sigma \d B(t)$ is added, trajectories
tend to saturate around those equilibrium points but jumping from one
equilibrium point to another is possible.  We use the MSE, the
mean$^2$+variance, the maximum-log-likelihood, and the proposed
finite-time-point time-decoupled squared $W_2$ distance
Eq.~\eqref{approximation} as loss functions to reconstruct
Eq.~\eqref{example1_ground_truth}.  For all loss functions, we use the
same neural network hyperparameters. Definitions of all loss functions
and training details are provided in Appendix~\ref{def_loss}. As
detailed in Appendix~\ref{training_details}, neural networks with the
same number of hidden layers and neurons in each layer are used for
each loss function. Using the initial condition $X(0)=0$, the sampled
ground-truth and reconstructed trajectories are shown in
Fig.~\ref{fig:fig1}.
\begin{figure}[h]
\centering 
\includegraphics[width=\textwidth]{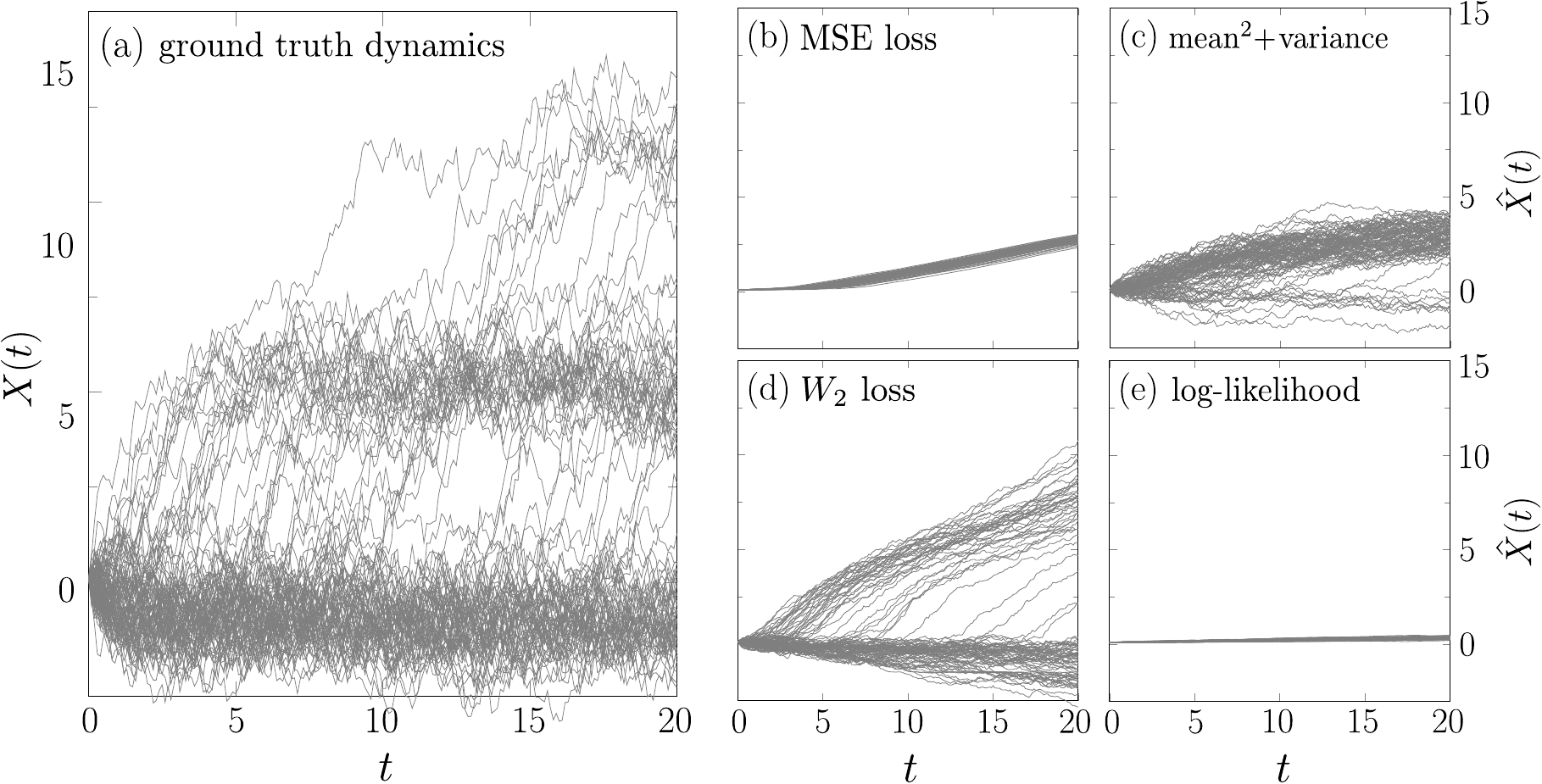}
\vspace{-0.2in}
\caption{(a) Ground-truth trajectories. (b) Reconstructed trajectories
  from nSDE using MSE loss. (c) Reconstructed trajectories from nSDE
  using mean$^2$+variance loss. (d) Reconstructed trajectories from
  nSDE using the finite-time-point time-decoupled $W_2$ loss.  (e)
  Reconstructed trajectories from nSDE using a max-log-likelihood loss
  yields the worst approximation.}
\label{fig:fig1}
\end{figure}

Fig.~\ref{fig:fig1}(a) shows the distributions of 100 trajectories
with most of them concentrated around two attractors (local minima
$x=-\tfrac{\pi}{3}, \tfrac{5\pi}{3}$ of the potential
$U(x)$). Fig.~\ref{fig:fig1}(b) shows that using MSE gives almost
deterministic trajectories and fails to reconstruct the
noise. From~\ref{fig:fig1}(c), we see that the mean$^2$+variance loss
fails to reconstruct the two local equilibria because cannot
sufficiently resolve the shape of the trajectory distribution at any
fixed timepoint. Fig.~\ref{fig:fig1}(d) shows that when using our
proposed finite-time-point time-decoupled squared $W_2$
loss~Eq.~\eqref{approximation}, the trajectories of the reconstructed
SDE can successfully learn the two-attractor feature and potentially
the distribution of trajectories. The reason why the reconstructed
trajectories of the $W_2$ distance cannot recover the third stable
equilibrium at $x=\tfrac{11\pi}{3}$ is because the data is sparse near
it. From~\ref{fig:fig1}(e), we see that the max-log-likelihood loss
performs the worst as it yields almost the same curves for all
realizations.

\end{example}

In the next example, we show how using our finite-time-point
time-decoupled squared $W_2$ distance loss function
Eq.~\eqref{approximation} can lead to efficient reconstruction of $f$
and $\sigma$. We shall use the mean relative $L^2$ error
\begin{equation}
  \Big(\sum_{i=0}^T\frac{\sum_{j=1}^N\|f(x_j(t_i), t_i)
    - \hat{f}(x_j(t_i), t_i)\|^2}{(T+1)\sum_{j=1}^N
    \|f(x_j(t_i), t_i)\|^2}\Big)^{\frac{1}{2}},\,\,
  \Big(\sum_{i=0}^T\frac{\sum_{j=1}^N\||\sigma(x_j(t_i), t_i)|
    - |\hat{\sigma}(x_j(t_i), t_i)|\|^2}{(T+1)\sum_{j=1}^N
    \|\sigma(x_j(t_i), t_i)\|^2}\Big)^{\frac{1}{2}}
\label{relative_l2}
\end{equation}
between the reconstructed $\hat{f}, \hat{\sigma}$ in
Eq.~\eqref{approximate_sde} and the ground-truth $f$ and $\sigma$ in
Eq.~\eqref{SDE_representation}. Here, $x_j(t_i)$ is the value of the
$j^{\text{th}}$ ground-truth trajectory at $t_i$.

\begin{example}
\rm
\label{example2}
Next, we reconstruct a Cox-Ingersoll-Ross (CIR) model which is a popular
finance model that describes the evolution of interest rates:
\begin{equation}
  \text{d} X(t) = \big(5 - X(t)\big)\text{d} t +
  \sigma_0\sqrt{X(t)}\text{d} B(t),\,\,\, t\in[0, 2].
\label{CIRmodel}
\end{equation}
Specifically, we are interested in how our reconstructed $\hat{f},
\hat{\sigma}$ can approximate the ground-truth $f(X) = 5 - X$ and
$\sigma(X) = \sigma_0\sqrt{X}$ (with $\sigma_0$ a constant
parameter). Here, we take the timestep $\Delta{t}=0.05$ in
Eq.~\eqref{approximation} and the initial condition is $X(0)=2$. For
reconstructing $f$ and $\sigma$, we compare using our proposed
finite-time-point time-decoupled squared $W_2$ distance
Eq.~\eqref{approximation} with minimizing a Maximum Mean Discrepancy
(MMD) \citep{briol2019statistical} and other loss functions given in
Appendix~\ref{def_loss}. Hyperparameters in the neural networks used
for training are the same across all loss functions.

\begin{figure}[h]
\centering 
\includegraphics[width=\textwidth]{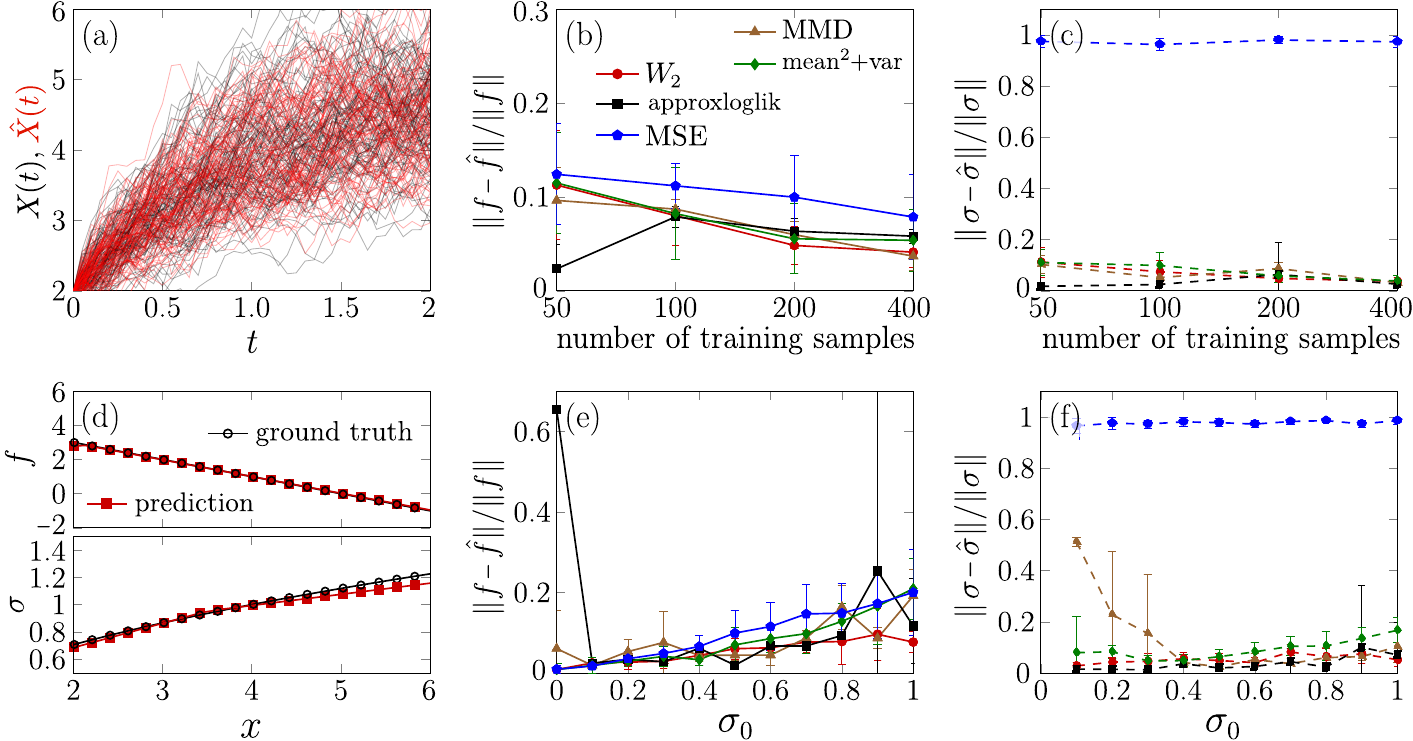}
\vspace{-0.2in}
\caption{(a) Ground-truth trajectories and reconstructed trajectories
  by nSDE using the finite-time-point time-decoupled squared $W_2$
  loss with $\sigma_0 = 0.5$. (b-c) Errors with respect to the numbers
  of ground-truth trajectories for $\sigma_0 = 0.5$. (d) Comparison of
  the reconstructed $\hat{f}_{\Theta_1}(u),
  \hat{\sigma}_{\Theta_2}(u)$ to the ground-truth functions $f(u),
  \sigma(u)$ for $\sigma_0 = 0.5$. (e-f) Errors with respect to noise
  level $\sigma_{0}$ with 200 training samples. Legends for panels (c,
  e, f) are the same as the one in (b).}
\label{fig:fig2}
\end{figure}
Fig.~\ref{fig:fig2}(a) shows the predicted trajectories using our
proposed squared $W_2$ loss function match well with the ground-truth
trajectories. Fig.~\ref{fig:fig2}(b, c) indicate that, if $\gtrsim
100$ ground-truth trajectories are used, our proposed squared $W_2$
distance loss yields smaller errors in $f, \sigma$ as defined in
Eq.~\eqref{relative_l2}. More specifically, we plot the reconstructed
$\hat{f}_{\Theta}, \hat{\sigma}_{\Theta}$ by using our squared $W_2$
loss in Fig.~\ref{fig:fig2}(d); these reconstructions also match well
with the ground-truth values $f, \sigma$. When we vary $\sigma_0$ in
Eq.~\eqref{CIRmodel}, our proposed finite-time-point time-decoupled
$W_2$ loss function gives the best performance among all loss
functions shown in Fig.~\ref{fig:fig2}(e, f). In
Appendix~\ref{IC_sensitivity}, instead of using the same initial
condition for all trajectories, we sample the initial condition from
different distributions and find that the reconstruction errors
$f-\hat{f}$ and $\sigma-\hat{\sigma}$ is \textbf{not} sensitive to
different initial conditions, implying the
robustness of using our proposed finite-time-point time-decoupled
$W_2$ loss function with respect to different initial
conditions. Also, in Appendix~\ref{neural_structure}, we change the
number of layers and the number of neurons in each layers for the two
neural networks we utilize to parameterize $\hat{f}\coloneqq \hat{f}(X,
t;\Theta_1), \hat{\sigma}\coloneqq \hat{\sigma}(X, t;\Theta_2)$. We
find that wider neural networks can lead to smaller errors
$f-\hat{f}$ and $\sigma-\hat{\sigma}$.

\end{example}

Next, we reconstruct the Ornstein-Uhlenbeck (OU) process given in
\cite{kidger2021neural} and in doing so, compare our loss function
with the WGAN-SDE method therein and with another recent MMD method.

\begin{example}
\rm
\label{example3}
Consider reconstructing the following time-inhomogeneous OU process
\begin{equation}
  \text{d} X(t) = \big(0.02t - 0.1 X(t)\big)\text{d} t
  + 0.4\text{d} B(t),\,\,\, t\in[0, 63].
\label{OUprocess}
\end{equation}
We compare the numerical performance of minimizing
Eq.~\eqref{time_discretize} or minimizing Eq.~\eqref{approximation}
with the WGAN method and using the MMD loss
metric. Eq.~\eqref{time_discretize} is numerically evaluated using the
\texttt{ot.emd2} function in the Python Optimal Transport package
\citep{flamary2021pot} We take the timestep $\Delta{t}=1$ in
Eq.~\eqref{approximation} and Eq.~\eqref{time_discretize} and the
initial condition is taken as $X(0)=0$. Neural networks with the same
number of hidden layers and neurons in each layer are used for all
three methods (see Appendix~\ref{training_details}).

In addition to the relative error in the reconstructed $\hat{f},
\hat{\sigma}$, we also compare the runtime and memory usage used by
the three methods as a function of the number of ground-truth
trajectories used in training.
\begin{figure}[h]
  \centering 
\includegraphics[width=\textwidth]{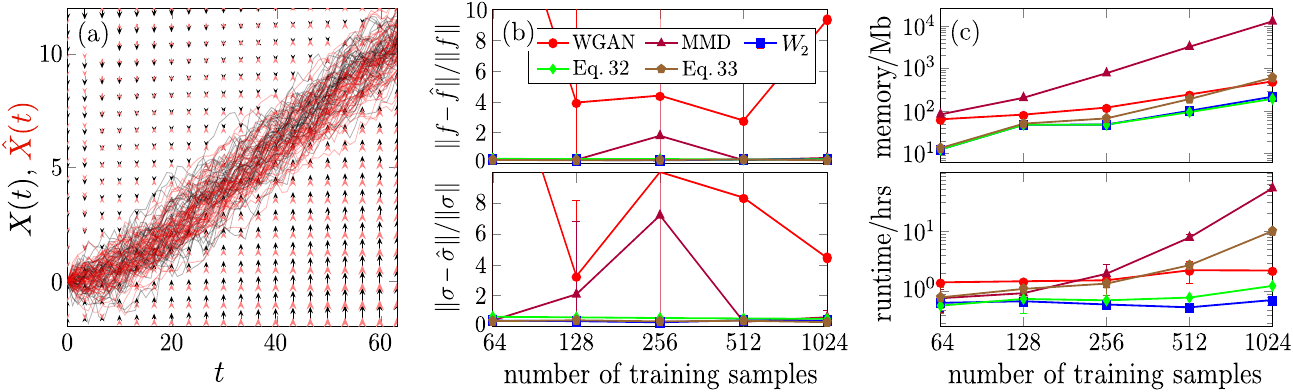}
\vspace{-0.22in}
\caption{(a) Ground-truth and reconstructed trajectories using the
  squared $W_2$ loss Eq.~\eqref{approximation}. Black and red curves
  are ground-truth and reconstructed trajectories, respectively. Black
  and red arrows indicate $f(x,t)$ and the reconstructed
  $\hat{f}(x,t)$ at fixed $(x,t)$, respectively. (b) Relative errors
  in reconstructed $\hat{f}$ and $\hat{\sigma}$, repeated 10
  times. Error bars show the standard deviation. (c) Resource
  consumption with respect to the number of training samples $N_{\rm
    samples}$. Memory usage is measured by torch profiler and
  represents peak memory usage during training. The legend in the
  panel (c) is the same as the one in (b).}
\label{fig:fig3}
\end{figure}
From Fig.~\ref{fig:fig3}(a), the distribution of trajectories of the
reconstructed SDE found from using our proposed squared $W_2$ loss
Eq.~\eqref{approximation} matches well with the distribution of the
ground-truth trajectories. Both minimizing Eq.~\eqref{time_discretize}
and minimizing Eq.~\eqref{approximation} outperform the other two
methods in the relative $L^2$ error of the reconstructed $f, \sigma$
for all numbers of ground-truth trajectories. Using
Eq.~\eqref{approximation} as the loss function achieves better
accuracy in a shorter computational time than using
Eq.~\eqref{time_discretize}.

For $N_{\text{sample}}$ training samples and $N$ total number of
timesteps, the memory cost in using Eq.~\eqref{approximation} is
$O(N\times N_{\text{sample}})$; however, the number of operations
needed is $O(N\times N_{\text{sample}}\log N_{\text{sample}})$ because
we need to reorder the ground-truth $X(t_i)$ and predicted
$\hat{X}(t_i)$ data to obtain the empirical cumulative distributions
at every $t_i$. The memory cost and operations needed in using
Eq.~\eqref{time_discretize} are both $O((N\times
N_{\text{sample}})^2)$ because a $(N\times N_{\text{sample}})\times
(N\times N_{\text{sample}})$ cost matrix must be evaluated.
On the other hand, the MMD method needs to create an
$N_{\text{sample}}\times N_{\text{sample}}$ matrix for each timestep
and thus the corresponding memory cost and operations needed are at
best $O(N\times N_{\text{sample}}^2)$. The WGAN-SDE method needs to
create a generator and a discriminator and its training is complex,
leading to both a higher memory cost and a larger runtime than our
method. For reconstructing SDEs, a larger number of ground-truth
trajectories leads to higher accuracy (see
Appendix~\ref{sgdtest}). Overall, our time-decoupled squared $W_2$
loss, Eq.~\eqref{approximation}, performs the best in terms of
accuracy and efficiency when reconstructing the 1D SDE
Eq.~\eqref{OUprocess}.

If we consider using stochastic gradient descent (SDG) to minibatch
for training, we find that the batch size cannot be set too small,
especially when we are using the MMD or Eq.~\eqref{time_discretize} as
loss functions, due to the intrinsic noisy nature of trajectories of
SDEs.  Thus, using our squared $W_2$ distance loss function given in
Eq.~\eqref{approximation} can be more efficient overall than using the
MMD or Eq.~\eqref{time_discretize} as the loss function. Additional
results using the SGD with minibatch for training are given in
Appendix~\ref{sgdtest}.

\end{example}

Finally, we carry out an experiment on reconstructing a 2D correlated
geometric Brownian motion. In this 2D reconstruction problem, we will
compare the loss functions, Eq.~\eqref{time_discretize} and
Eq.~\eqref{approximation}, the MMD method, and a sliced squared
Wasserstein distance method \citep{kolouri2018sliced}.

\begin{example}
\rm
\label{example4}
Consider reconstructing the following 2D correlated geometric
Brownian motion that can represent, \textit{e.g.}, values of two correlated
stocks \citep{musiela2006martingale}

\begin{equation}
  \begin{aligned}
    \text{d} X_1(t) = & \mu_1 X_1(t)\d t + \sum_{i=1}^2 \sigma_{1, i}X_i(t)\d B_i(t),\\
    \text{d} X_2(t) = & \mu_2 X_2(t)\d t + \sum_{i=1}^2 \sigma_{2, i}X_i(t)\d B_i(t)
    \end{aligned}
\label{correlated_Brown}
\end{equation}
Here, $t\in[0, 2]$, $B_1(t)$ and $B_2(t)$ are independent Brownian
processes, $\bm{f} \coloneqq (\mu_1 X_1, \mu_2 X_2)$ is a 2D vector,
and $\bm{\sigma}\coloneqq[\sigma_{1, 1}X_1, \sigma_{1, 2}X_2;
  \sigma_{2, 1}X_1, \sigma_{2, 2}X_2]$ is a $2\times 2$ matrix.
\begin{figure}[h]
\includegraphics[width=0.75\textwidth]{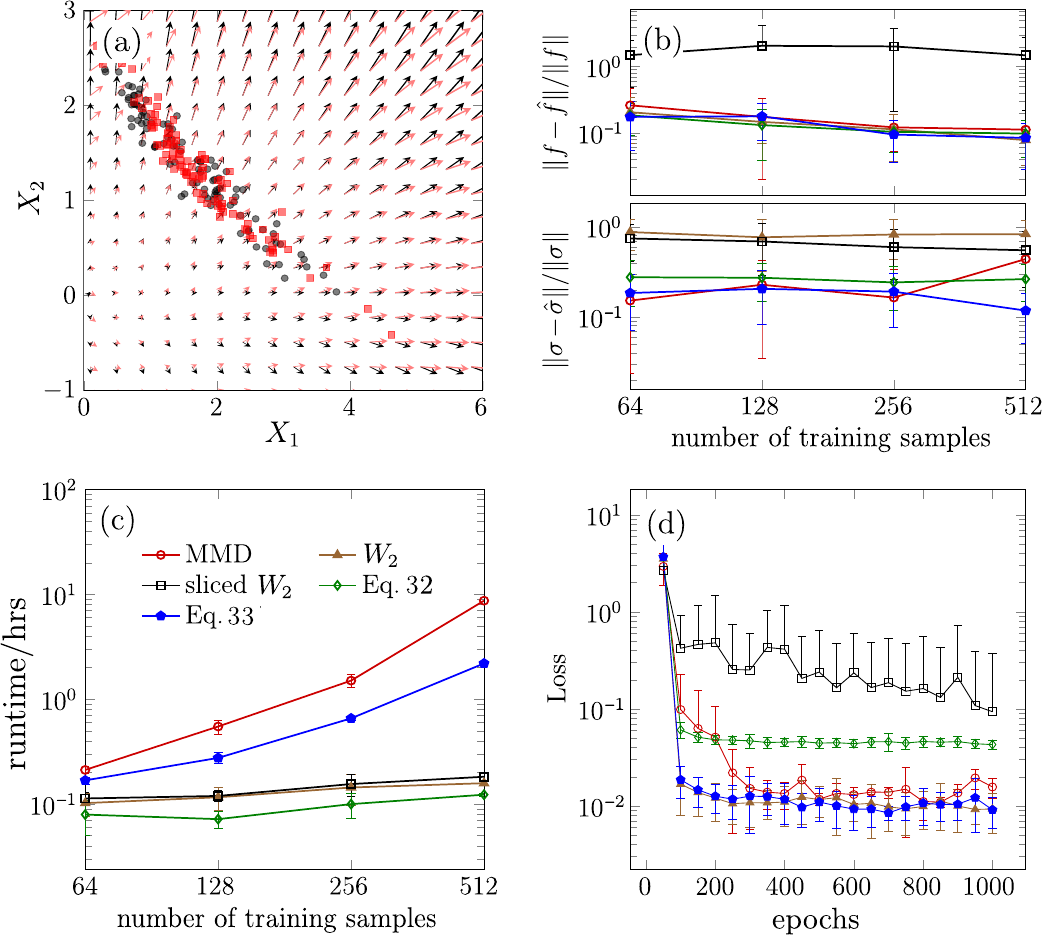}
\centering 
\caption{(a) Black dots and red squares are the ground-truth
  $(X_1(2), X_2(2))$ and the reconstructed $(\hat{X}_1(2),
  \hat{X}_2(2))$ found using the rotated squared $W_2$ loss function,
  respectively. Black and red arrows indicate, respectively, the
  vectors ${\bm f}(X_1,X_2)$ and $\hat{{\bm f}}(X_1,X_2)$. (b)
  Relative errors of the reconstructed ${\bm f}$ and
  $\boldsymbol{\sigma}$. Error bars indicate the standard deviation
  across ten reconstructions. (c) Runtime of different loss functions
  with respect to $N_{\rm samples}$. (d) The decrease of different
  loss functions with respect to training epochs. The legend for the
  panel (d) is the same as the one in (c).}
\vspace{-8mm}
\label{fig:fig4}
\end{figure}
We use $(\mu_1, \mu_2) = (0.1, 0.2)$, $\bm{\sigma}=[0.2X_1, -0.1X_2;
  -0.1X_1, 0.1X_2]$, and set the initial condition $(X_1(0), X_2(0)) =
(1, 0.5)$. In addition to directly minimizing a 2D decorrelated
version of the squared $W_2$ distance Eq.~\eqref{approximation}
(denoted as $W_2$ in Fig.~\ref{fig:fig4}(c)), we consider minimizing a
sliced squared $W_2$ distance as proposed by
\cite{kolouri2018sliced,kolouri2019generalized}.  Finally, we
numerically estimate the $W_2$ distance Eq.~\eqref{time_discretize} as
well as the time-decoupled approximation Eq.~\eqref{approximation}
using the \texttt{ot.emd2} function in the Python Optimal Transport
package.  Formulae of the above loss functions are given in
Appendix~\ref{def_loss}. We keep the neural network hyperparameters
the same while minimizing all loss functions.  Note that since the SDE
has two components, the definition of the relative error in $\sigma$
is revised to
\begin{equation}
\left[\sum_{i=0}^{T} \frac{\sum_{j=1}^{N}
\|\bm{\sigma}\bm{\sigma}^T(x_j(t_i), t_i) -
\hat{\bm{\sigma}}\hat{\bm{\sigma}}^T(x_j(t_i),
t_i)\|^2_F}{(T+1)\sum_{j=1}^N\|\hat{\bm{\sigma}}\hat{\bm{\sigma}}^T(x_{j}(t_i),
t_i)\|_F^2}\right]^{1/2},
\label{sigma_error}
\end{equation}
where $\|\cdot\|_F$ is the Frobenius norm for matrices.

Fig.~\ref{fig:fig4}(a) shows the ground truth and reconstructed
coordinates $(X_1, X_2)$ (black dots) and $(\hat{X}_1, \hat{X}_2)$
(red squares) at time $t=2$, along with ${\bm f}(X_1,X_2)$ (black) and
$\hat{{\bm f}}(X_1,X_2)$ (red). For reconstructing ${\bm f}$ and
$\bm{\sigma}$ in problem, numerically evaluating
Eq.~\eqref{approximation} (blue curve) performs better than the MMD
method, the loss in Eq.~\eqref{time_discretize}, the sliced $W_2$
distance loss, and the 2D decorrelated squared $W_2$ loss, as shown in
Fig.~\ref{fig:fig4}(b). Using the sliced $W_2$ distance yields the
poorest performance and least accurate $\hat{f}$ and
$\hat{\bm{\sigma}}$. Using the 2D decorrelated squared $W_2$ loss
function also gives inaccurate $\hat{\bm{\sigma}}$. Thus, the sliced
$W_2$ distance and the 2D decorrelated squared $W_2$ loss are not good
candidates for reconstructing multivariate SDEs.  Numerically
estimating Eq.~\eqref{time_discretize} yields poorer performance than
numerically estimating Eq.~\eqref{approximation} because numerically
evaluating the $W_2$ distance for higher-dimensional empirical
distributions is generally less accurate.

From Fig.~\ref{fig:fig4}(c), we see that the runtime and memory needed
to numerically evaluate the time-decoupled~Eq.~\eqref{approximation}
using \texttt{ot.emd2} is smaller than those needed for the MMD
method, but larger than those needed to numerically estimate
Eq.~\eqref{time_discretize}. Yet, as shown in Fig.~\ref{fig:fig4}(d),
minimizing Eq.~\eqref{approximation} leads to the fastest convergence,
potentially requiring fewer epochs when using Eq.~\eqref{approximation} as
the loss function.  An additional comparison of using the two loss functions, the finite-time-point squared $W_2$ distance Eq.~\eqref{time_discretize}
and the finite-time-point time-decoupled squared $W_2$ distance Eq.~\eqref{approximation}
is given in Appendix~\ref{appendix_more_discussion}. Further analysis on how the number of samples and
the dimensionality of an SDE dimensionality affects $W_2$-based
distances in reconstructing multivariate SDEs will be informative.

\end{example}
\section{Summary and Conclusions}
\label{summary}
In this paper, we analyzed the squared $W_2$ distance between two
probability distributions associated with two SDEs and proposed a
novel method for efficient reconstruction of SDEs from data by
minimizing squared $W_2$ distances as loss functions.  Upon performing
numerical experiments, we found that our proposed finite-time-point
time-decoupled squared $W_2$ distance loss function,
Eq.~\eqref{approximation}, is superior than many other recently
developed machine-learning and statistical approaches to SDE
reconstruction.

A number of extensions are apparent. First, one can further
investigate applying the squared $W_2$ loss to the reconstruction of
high-dimensional SDEs.  Another important direction to develop are
approaches to efficiently evaluate the squared $W_2$ loss function
Eq.~\eqref{time_discretize} and analyze how well the time-decoupled
squared $W_2$ loss function Eq.~\eqref{approximation} can approximate
Eq.~\eqref{time_discretize}. Whether the Wasserstein distance can
serve as upper bounds for the errors $f-\hat{f}$ and
$\sigma-\hat{\sigma}$ is also an intriguing question as its resolution
will determine is minimizing the squared Wasserstein distance is
sufficient for reconstructing SDEs.  Finally, another promising area
worthy of study is the extension of the squared $W_2$ distance loss
function to the reconstruction of general L\'evy processes that
include jumps in the trajectories.


\acks{The authors acknowledge Dr. Hedi Xia for his constructive
  comments on this manuscript. TC acknowledges support from the Army
  Research Office through grant W911NF-18-1-0345.}


\vskip 0.2in
\bibliographystyle{plain}
\bibliography{iclr2024_conference}

\begin{thebibliography}{41}
\providecommand{\natexlab}[1]{#1}
\providecommand{\url}[1]{\texttt{#1}}
\expandafter\ifx\csname urlstyle\endcsname\relax
  \providecommand{\doi}[1]{doi: #1}\else
  \providecommand{\doi}{doi: \begingroup \urlstyle{rm}\Url}\fi

\bibitem[Arjovsky et~al.(2017)Arjovsky, Chintala, and
  Bottou]{arjovsky2017wasserstein}
Martin Arjovsky, Soumith Chintala, and L{\'e}on Bottou.
\newblock {W}asserstein generative adversarial networks.
\newblock In \emph{International Conference on Machine Learning}, pages
  214--223. PMLR, 2017.

\bibitem[Bartl et~al.(2021)Bartl, Beiglb{\"o}ck, and
  Pammer]{bartl2021wasserstein}
Daniel Bartl, Mathias Beiglb{\"o}ck, and Gudmund Pammer.
\newblock The {W}asserstein space of stochastic processes.
\newblock \emph{arXiv preprint arXiv:2104.14245}, 2021.

\bibitem[Bion-Nadal and Talay(2019)]{bion2019wasserstein}
Jocelyne Bion-Nadal and Denis Talay.
\newblock On a {W}asserstein-type distance between solutions to stochastic
  differential equations.
\newblock \emph{The Annals of Applied Probability}, 29\penalty0 (3):\penalty0
  1609--1639, 2019.
\newblock \doi{10.1214/18-aap1423}.

\bibitem[Bressloff(2014)]{bressloff2014stochastic}
Paul~C Bressloff.
\newblock \emph{Stochastic Processes in Cell Biology}, volume~41.
\newblock Springer, 2014.

\bibitem[Briol et~al.(2019)Briol, Barp, Duncan, and
  Girolami]{briol2019statistical}
Francois-Xavier Briol, Alessandro Barp, Andrew~B Duncan, and Mark Girolami.
\newblock Statistical inference for generative models with maximum mean
  discrepancy.
\newblock \emph{arXiv preprint arXiv:1906.05944}, 2019.

\bibitem[Chen et~al.(2018)Chen, Rubanova, Bettencourt, and
  Duvenaud]{chen2018neural}
Ricky~TQ Chen, Yulia Rubanova, Jesse Bettencourt, and David~K Duvenaud.
\newblock Neural ordinary differential equations.
\newblock In \emph{Advances in Neural Information Processing Systems},
  volume~31, 2018.

\bibitem[Chewi et~al.(2021)Chewi, Clancy, Le~Gouic, Rigollet, Stepaniants, and
  Stromme]{chewi2021fast}
Sinho Chewi, Julien Clancy, Thibaut Le~Gouic, Philippe Rigollet, George
  Stepaniants, and Austin Stromme.
\newblock Fast and smooth interpolation on {W}asserstein space.
\newblock In \emph{International Conference on Artificial Intelligence and
  Statistics}, pages 3061--3069. PMLR, 2021.

\bibitem[Cinlar(2011)]{cinlar2011probability}
Erhan Cinlar.
\newblock \emph{Probability and Stochastics}.
\newblock Springer Science \& Business Media, 2011.

\bibitem[Clement and Desch(2008)]{clement2008elementary}
Philippe Clement and Wolfgang Desch.
\newblock An elementary proof of the triangle inequality for the {W}asserstein
  metric.
\newblock \emph{Proceedings of the American Mathematical Society}, 136\penalty0
  (1):\penalty0 333--339, 2008.

\bibitem[Cuturi et~al.(2019)Cuturi, Teboul, and Vert]{cuturi2019differentiable}
Marco Cuturi, Olivier Teboul, and Jean-Philippe Vert.
\newblock Differentiable ranks and sorting using optimal transport.
\newblock In \emph{Proceedings of the 33rd International Conference on Neural
  Information Processing Systems}, pages 6861--6871, 2019.

\bibitem[De~Vecchi et~al.(2016)De~Vecchi, Morando, and
  Ugolini]{de2016reduction}
Francesco~C De~Vecchi, Paola Morando, and Stefania Ugolini.
\newblock Reduction and reconstruction of stochastic differential equations via
  symmetries.
\newblock \emph{Journal of Mathematical Physics}, 57\penalty0 (12), 2016.

\bibitem[Flamary et~al.(2021)Flamary, Courty, Gramfort, Alaya, Boisbunon,
  Chambon, Chapel, Corenflos, Fatras, Fournier, et~al.]{flamary2021pot}
R{\'e}mi Flamary, Nicolas Courty, Alexandre Gramfort, Mokhtar~Z Alaya,
  Aur{\'e}lie Boisbunon, Stanislas Chambon, Laetitia Chapel, Adrien Corenflos,
  Kilian Fatras, Nemo Fournier, et~al.
\newblock {Pot: Python optimal transport}.
\newblock \emph{The Journal of Machine Learning Research}, 22\penalty0
  (1):\penalty0 3571--3578, 2021.

\bibitem[Fournier and Guillin(2015)]{fournier2015rate}
Nicolas Fournier and Arnaud Guillin.
\newblock On the rate of convergence in {W}asserstein distance of the empirical
  measure.
\newblock \emph{Probability Theory and Related Fields}, 162\penalty0
  (3-4):\penalty0 707--738, 2015.

\bibitem[Frogner et~al.(2015)Frogner, Zhang, Mobahi, Araya, and
  Poggio]{frogner2015learning}
Charlie Frogner, Chiyuan Zhang, Hossein Mobahi, Mauricio Araya, and Tomaso~A
  Poggio.
\newblock Learning with a {W}asserstein loss.
\newblock In \emph{Advances in Neural Information Pprocessing Systems},
  volume~28, 2015.

\bibitem[He et~al.(2016)He, Zhang, Ren, and Sun]{he2016deep}
Kaiming He, Xiangyu Zhang, Shaoqing Ren, and Jian Sun.
\newblock Deep residual learning for image recognition.
\newblock In \emph{Proceedings of the IEEE conference on computer vision and
  pattern recognition}, pages 770--778, 2016.

\bibitem[Jia and Benson(2019)]{jia2019neural}
Junteng Jia and Austin~R Benson.
\newblock Neural jump stochastic differential equations.
\newblock \emph{Advances in Neural Information Processing Systems}, 32, 2019.

\bibitem[Kidger et~al.(2021)Kidger, Foster, Li, and Lyons]{kidger2021neural}
Patrick Kidger, James Foster, Xuechen Li, and Terry~J Lyons.
\newblock Neural {SDE}s as infinite-dimensional {GAN}s.
\newblock In \emph{International Conference on Machine Learning}, pages
  5453--5463. PMLR, 2021.

\bibitem[Kolouri et~al.(2018)Kolouri, Rohde, and Hoffmann]{kolouri2018sliced}
Soheil Kolouri, Gustavo~K Rohde, and Heiko Hoffmann.
\newblock Sliced {W}asserstein distance for learning {G}aussian mixture models.
\newblock In \emph{Proceedings of the IEEE Conference on Computer Vision and
  Pattern Recognition}, pages 3427--3436, 2018.

\bibitem[Kolouri et~al.(2019)Kolouri, Nadjahi, Simsekli, Badeau, and
  Rohde]{kolouri2019generalized}
Soheil Kolouri, Kimia Nadjahi, Umut Simsekli, Roland Badeau, and Gustavo Rohde.
\newblock Generalized sliced {W}asserstein distances.
\newblock In \emph{Advances in Neural Information Processing Systems},
  volume~32, 2019.

\bibitem[Li et~al.(2020)Li, Wong, Chen, and Duvenaud]{li2020scalable}
Xuechen Li, Ting-Kam~Leonard Wong, Ricky~TQ Chen, and David Duvenaud.
\newblock Scalable gradients for stochastic differential equations.
\newblock In \emph{International Conference on Artificial Intelligence and
  Statistics}, pages 3870--3882. PMLR, 2020.

\bibitem[Li et~al.(2015)Li, Swersky, and Zemel]{li2015generative}
Yujia Li, Kevin Swersky, and Rich Zemel.
\newblock Generative moment matching networks.
\newblock In \emph{International Conference on Machine Learning}, pages
  1718--1727. PMLR, 2015.

\bibitem[Lin and Buchler(2018)]{lin2018efficient}
Yen~Ting Lin and Nicolas~E Buchler.
\newblock Efficient analysis of stochastic gene dynamics in the non-adiabatic
  regime using piecewise deterministic {M}arkov processes.
\newblock \emph{Journal of The Royal Society Interface}, 15\penalty0
  (138):\penalty0 20170804, 2018.

\bibitem[Liu et~al.(2020)Liu, Ong, Shen, and Cai]{liu2020gaussian}
Haitao Liu, Yew-Soon Ong, Xiaobo Shen, and Jianfei Cai.
\newblock When {G}aussian process meets big data: a review of scalable {GP}s.
\newblock \emph{IEEE Transactions on Neural Networks and Learning Systems},
  31\penalty0 (11):\penalty0 4405--4423, 2020.

\bibitem[MacKay et~al.(1998)]{mackay1998introduction}
David~JC MacKay et~al.
\newblock Introduction to {G}aussian processes.
\newblock \emph{NATO ASI Series F Computer and Systems Sciences}, 168:\penalty0
  133--166, 1998.

\bibitem[Mathelin et~al.(2005)Mathelin, Hussaini, and
  Zang]{mathelin2005stochastic}
Lionel Mathelin, M~Yousuff Hussaini, and Thomas~A Zang.
\newblock Stochastic approaches to uncertainty quantification in {CFD}
  simulations.
\newblock \emph{Numerical Algorithms}, 38:\penalty0 209--236, 2005.

\bibitem[Musiela and Rutkowski(2006)]{musiela2006martingale}
Marek Musiela and Marek Rutkowski.
\newblock \emph{Martingale Methods in Financial Modelling}, volume~36.
\newblock Springer Science \& Business Media, 2006.

\bibitem[Oh et~al.(2019)Oh, Pouryahya, Iyer, Apte, Tannenbaum, and
  Deasy]{oh2019kernel}
Jung~Hun Oh, Maryam Pouryahya, Aditi Iyer, Aditya~P Apte, Allen Tannenbaum, and
  Joseph~O Deasy.
\newblock Kernel {W}asserstein distance.
\newblock \emph{arXiv preprint arXiv:1905.09314}, 2019.

\bibitem[Pereira et~al.(2010)Pereira, Ibrahimi, and
  Montanari]{pereira2010learning}
Jos{\'e} Pereira, Morteza Ibrahimi, and Andrea Montanari.
\newblock Learning networks of stochastic differential equations.
\newblock In \emph{Advances in Neural Information Processing Systems},
  volume~23, 2010.

\bibitem[Rowland et~al.(2019)Rowland, Hron, Tang, Choromanski, Sarlos, and
  Weller]{rowland2019orthogonal}
Mark Rowland, Jiri Hron, Yunhao Tang, Krzysztof Choromanski, Tamas Sarlos, and
  Adrian Weller.
\newblock Orthogonal estimation of {W}asserstein distances.
\newblock In \emph{The 22nd International Conference on Artificial Intelligence
  and Statistics}, pages 186--195. PMLR, 2019.

\bibitem[R{\"u}schendorf(1985)]{ruschendorf1985wasserstein}
Ludger R{\"u}schendorf.
\newblock The {W}asserstein distance and approximation theorems.
\newblock \emph{Probability Theory and Related Fields}, 70\penalty0
  (1):\penalty0 117--129, 1985.

\bibitem[Sanz-Serna and Zygalakis(2021)]{sanz2021wasserstein}
Jesus~Maria Sanz-Serna and Konstantinos~C Zygalakis.
\newblock {W}asserstein distance estimates for the distributions of numerical
  approximations to ergodic stochastic differential equations.
\newblock \emph{The Journal of Machine Learning Research}, 22\penalty0
  (1):\penalty0 11006--11042, 2021.

\bibitem[Soize(2017)]{soize2017uncertainty}
Christian Soize.
\newblock \emph{Uncertainty Quantification}.
\newblock Springer, 2017.

\bibitem[Song et~al.(2020)Song, Sohl-Dickstein, Kingma, Kumar, Ermon, and
  Poole]{song2020score}
Yang Song, Jascha Sohl-Dickstein, Diederik~P Kingma, Abhishek Kumar, Stefano
  Ermon, and Ben Poole.
\newblock Score-based generative modeling through stochastic differential
  equations.
\newblock In \emph{International Conference on Learning Representations}, 2020.

\bibitem[Sullivan(2015)]{sullivan2015introduction}
Timothy~John Sullivan.
\newblock \emph{Introduction to Uncertainty Quantification}, volume~63.
\newblock Springer, 2015.

\bibitem[Tong et~al.(2022)Tong, Nguyen-Tang, Tran, and Choi]{tong2022learning}
Anh Tong, Thanh Nguyen-Tang, Toan Tran, and Jaesik Choi.
\newblock Learning fractional white noises in neural stochastic differential
  equations.
\newblock In \emph{Advances in Neural Information Processing Systems},
  volume~35, pages 37660--37675, 2022.

\bibitem[Tzen and Raginsky(2019)]{tzen2019neural}
Belinda Tzen and Maxim Raginsky.
\newblock Neural stochastic differential equations: deep latent {G}aussian
  models in the diffusion limit.
\newblock \emph{arXiv preprint arXiv:1905.09883}, 2019.

\bibitem[Villani et~al.(2009)]{villani2009optimal}
C{\'e}dric Villani et~al.
\newblock \emph{Optimal transport: old and new}, volume 338.
\newblock Springer, 2009.

\bibitem[Wang(2016)]{wang2016p}
Jian Wang.
\newblock ${L}^{p}$-{W}asserstein distance for stochastic differential
  equations driven by {L\'e}vy processes.
\newblock \emph{Bernoulli}, pages 1598--1616, 2016.

\bibitem[Welch et~al.(1995)Welch, Bishop, et~al.]{welch1995introduction}
Greg Welch, Gary Bishop, et~al.
\newblock An introduction to the {K}alman filter.
\newblock 1995.

\bibitem[Welch(2020)]{welch2020kalman}
Gregory~F Welch.
\newblock Kalman filter.
\newblock \emph{Computer Vision: A Reference Guide}, pages 1--3, 2020.

\bibitem[Zheng et~al.(2020)Zheng, Wang, and Gou]{zheng2020nonparametric}
Wenbo Zheng, Fei-Yue Wang, and Chao Gou.
\newblock Nonparametric different-feature selection using {W}asserstein
  distance.
\newblock In \emph{2020 IEEE 32nd International Conference on Tools with
  Artificial Intelligence (ICTAI)}, pages 982--988. IEEE, 2020.

\end{thebibliography}

\newpage

\appendix
\section{Proof to Theorem~\ref{theorem1}}
\label{proof_theorem1}
Here, we shall provide a proof to Theorem~\ref{theorem1}.
First, note that $\tilde{X}(t)$ defined in Eq.~\eqref{tilde_x} is a
specific realization of $\hat{X}(t)$ defined in
Eq.~\eqref{approximate_sde}, \textit{coupled to} $X(t)$ in the sense
that its initial values are $X(0)$ almost surely and the It{\^o}
integral is defined with respect to the same standard Brownian motion
$B(t)$. Therefore, by definition, if we let $\pi$ in Eq.~\eqref{pidef}
to be the joint distribution of $(X, \tilde{X})$, then
\begin{equation}
  W_2(\mu, \hat \mu) \leq \left(\E\Big[\medint\int_0^T|\tilde{X}(t)
    - X(t)|^2\d t\Big]\right)^{1/2}.
\end{equation}
Next, we provide a bound for $\E\big[\int_0^T|\tilde{X}(t) - X(t)|^2\d
  t\big]^{\frac{1}{2}}$ by the mean value theorem for $f$ and $g$.
\begin{equation}
\begin{aligned}
\d \big(X(t)-\tilde{X}(t)\big) = &\,  \partial_x f\big(\eta_1(X(t), \tilde{X}(t),t),
  t\big)\cdot (X(t)-\tilde{X}(t))\d t \\
\: & + \partial_x\sigma\big(\eta_2 (X(t), \tilde{X}(t)), t\big)
\cdot (X(t)-\tilde{X}(t))\d B(t) \\
\: & + (f-\hat{f})(\tilde{X}(t), t)\d t
+ (\sigma-\hat{\sigma})(\tilde{X}(t), t))\d B(t). 
\end{aligned}
\end{equation}
where $\eta_1(x_1, x_2), \eta_{2}(x_1, x_2)$ are defined in
Eq.~\eqref{intermediate_eq} such that their values are in $(x_1,
x_2)$.

Applying It$\hat{\rm{o}}$'s formula to $[X(t) -\tilde X(t)]/H(0;t)$,
we find
\begin{equation}
\begin{aligned}
\text{d}\left(\frac{X(t)-\tilde{X}(t)}{H(0;t)}\right)
= & \frac{1}{H(0;t)}
\Big[(f-\hat{f})(\tilde{X}(t),t)\d t +\partial_x\sigma\big(\eta_2 (X,\tilde{X}), t\big)\cdot
  (\sigma-\hat{\sigma})(\tilde{X}(t), t) \d t\Big]\\
\: & \qquad + \frac{1}{H(0;t)}\Big[(\sigma-\hat{\sigma})(\tilde{X}(t), t)\d B(t)\Big].
\end{aligned}
\end{equation}
Integrating both sides from $0$ to $t$, we obtain
\begin{equation}
\begin{aligned}
X(t) - \tilde{X}(t) 
= & \int_0^t H(s;t) \Big[(f-\hat{f})(\tilde{X}(s),s) +
  \partial_x\sigma\big(\eta_2 (X,\tilde{X}), s\big)
  \cdot (\sigma-\hat{\sigma})(\tilde{X}(s), s)\Big]\d s \\
\: & \qquad +\int_0^t  H(s;t)
\cdot(\sigma-\hat{\sigma})(\tilde{X}(s), s)\d B(s).
\end{aligned}
\end{equation}
By invoking It$\hat{\rm{o}}$ isometry and observing that $(a+b+c)^2 \leq
3(a^2 + b^2 + c^2)$, we deduce

\begin{equation}
\begin{aligned}
& \E\big[\big(X(t) - \tilde{X}(t)\big)^2\big] \leq 3 \E\Big[\Big(\medint\int_0^t H(s;t)
\cdot (f-\hat{f})(\tilde{X}(s),s)\d s\Big)^2\Big] \\
&\quad\quad\quad\quad+3 \E\Big[\Big(\medint\int_0^t H(s;t)
    \cdot (\partial_x\sigma\big(\eta_2 (X,\tilde{X}), s\big)\cdot
    (\sigma-\hat{\sigma})(\tilde{X}(s), s)\d s \Big)^2\Big] \\
  &\quad\quad\quad\quad
  +3\E\Big[\Big(\medint\int_0^t H(s;t)\cdot(\sigma-\hat{\sigma})
    (\tilde{X}(s), s)\d B(s)\Big)^2\Big]\\
  &\quad\quad\leq 3\E\Big[\medint\int_0^t H^{2}(s;t)\d s\Big]
  \times \E\Big[\medint\int_0^T (f- \hat{f})^2(\tilde{X}(s), s)\d s\Big] \\ 
&\quad\quad\quad\quad+ 3\E\Big[\medint\int_0^t H^{2}(s;t)\d s\Big] \times
  \E\Big[\medint\int_0^T \Big(\partial_x\sigma\big(\eta_2 (X,\tilde{X}), s\big)
    \cdot (\sigma-\hat{\sigma})(\tilde{X}(s), s)\Big)^2\d s\Big]\\
  &\quad\quad\quad\quad + 3\E\Big[\medint\int_0^t H^{2}(s;t)
    \cdot (\sigma - \tilde{\sigma})^2(\tilde{X}(s), s)\d s\Big]\\
  &\quad\quad\leq 3\E\Big[\medint\int_0^t H^{2}(s;t) \d s\Big]\times
  \E\Big[\medint\int_0^t (f- \hat{f})^2(\tilde{X}(s), s)\d s\Big] \\
&\quad\quad\quad\quad+3\E\Big[\medint\int_0^t H^{2}(s;t)  \d s\Big] \times 
  \E\Big[\medint\int_0^t \Big(\partial_x\sigma\big(\eta_2 (X,\tilde{X}), s\big)
    \cdot (\sigma-\hat{\sigma})(\tilde{X}(s), s)\Big)^2\text{d}s \Big]\\
  &\quad\quad\quad\quad+
  3\left(\E\Big[\medint\int_0^t H^{4}(s;t) \d s\Big]\right)^{1/2}\!\!
  \times \left(\E\Big[\medint\int_0^t
    (\sigma- \hat{\sigma})^4(\tilde{X}(s), s)\d s\Big]\right)^{1/2}.
\end{aligned}
\end{equation}
Finally, we conclude that

\begin{equation}
\begin{aligned}
  W_2^2(\mu, \tilde{\mu}) \leq & \int_0^T
  \E\big[\big(X(t)-\tilde{X}(t)\big)^2\big]\d t\\
\leq & 3\int_0^T \E\Big[\medint\int_0^t H^{2}(s;t) \d s\Big]\d t
\times \E\Big[\medint\int_0^T (f- \hat{f})^2(\tilde{X}(s), s)\d s)\Big]\\
\: & \, +3\int_0^T\E\Big[\medint\int_0^t H^{2}(s;t) \d s\Big]\d t \times 
\E\Big[\medint\int_0^T (\partial_x\sigma(\eta_2 (X,\tilde{X}), s)\cdot
  (\sigma-\hat{\sigma})(\tilde{X}(s), s)\large)^2 \d s)\Big]\\
\: &\, +3\int_0^T \left(\E\Big[\medint\int_0^t H^{4}(s;t) \d s\Big]\right)^{1/2}\d t
\times \left(\E\Big[\medint\int_0^T
  (\sigma- \hat{\sigma})^4(\tilde{X}(s), s)\d s)\Big]\right)^{1/2},
\end{aligned}
\end{equation}
which proves Theorem~\ref{theorem1}.

This theorem can be generalized to higher dimensional dynamics of
$\bm{X}(t)$ and $\hat{\bm{X}}(t)$ described by

\begin{equation}
\begin{aligned}
\text{d} \bm{X}(t) = & \bm{f}(\bm{X}(t), t)\text{d}t + 
\bm{\sigma}(\bm{X}(t), t)\text{d} \bm{B}(t),\\
\text{d} \hat{\bm{X}}(t) = & \hat{\bm{f}}(\hat{\bm{X}}(t), t)\text{d}t 
+ \hat{\bm{\sigma}}(\hat{\bm{X}}(t), t)\text{d} \hat{\bm{B}}(t),
\end{aligned}
\label{multi_dimensional}
\end{equation}
where $\bm{f},\hat{\bm{f}}: \mathbb{R}^{d+1}\rightarrow\mathbb{R}^d$
are the $d$-dimensional drift functions and
$\bm{\sigma},\hat{\bm{\sigma}}:
\mathbb{R}^{d+1}\rightarrow\mathbb{R}^{d\times {s}}$ are diffusion
matrices. $\bm{B}(t)$ and $\bm{\hat{B}}(t)$ are two independent
${s}$-dimensional standard Brownian motions. Under some additional
assumptions, one can derive an upper bound for the $W_2$ distance of
the probability distributions for $\bm{X}, \bm{\hat{X}}$ as in
Eq.~\eqref{W2_estimate}. For example, if for every $i=1,...d$, the
$i^{\text{th}}$ component $\d X_{i}(t) = f_{i}(X_{i}(t), t)\d t +
\sigma_{i}(X_{i}(t), t) \d B_{i}(t)$ and $\text{d} \hat{X}_{i}(t) =
\hat{f}_{i}(\hat{X}_{i}(t), t)\d t + \hat{\sigma}_{i}(\hat{X}_{i}(t),
t) \text{d} \hat{B}_{i}(t)$, where $B_i(t), \hat{B}_i(t)$ are
independent Brownian motions, then similar conclusions can be derived
by calculating the difference $X_{i}-\hat{X}_{i}$. Developing an upper
bound of the $W_2$ distance for general dimensions $d$ requires
additional assumptions to find expressions for
$\bm{X}-\hat{\bm{X}}$. We leave this nontrivial derivation as future
work. Although without a formal theoretical analysis, we shall show in
Example \ref{example4} that applying the squared $W_2$ distance as the
loss function is also effective in reconstructing multidimensional
SDEs.

\section{Single-trajectory MSE and KL divergence}
\label{proof_theorem2}
We shall first show that using the single-trajectory MSE tends to fit
the mean process $\E[X(t)]$ and make noise diminish, which indicates
that the MSE is not a good loss function when one wishes to fit
$\sigma$ in Eq.~\eqref{SDE_representation}.

For two \textit{independent} $d$-dimensional stochastic processes
$\{\bm{X}(t)\}_{t=0}^T, \{\hat{\bm{X}}(t)\}_{t=0}^T$ as solutions to
Eq.~\eqref{multi_dimensional} with appropriate ${\bm f}, \hat{\bm f}$
and ${\bm \sigma}, \hat{\bm \sigma}$, let $\E[{\bm X}]$ represent the
trajectory of mean values of $\bm{X}(t)$, \textit{i.e.}, $\E[\bm{X}] =
\E[\bm{X}(t)]$. We have
\begin{equation}
\begin{aligned}
  \E \big[\| \boldsymbol{X} - \hat {\boldsymbol{X}}\|^2\big]
  = & \E\big[\left\|\boldsymbol{X}- \E[\boldsymbol{X}]\right\|^2\big]
  + \E\big[\|\hat{\boldsymbol{X}} - \E[\boldsymbol{X}]\|^2\big] \\
  \: & \quad - 2\E\Big[\medint\int_0^T \!\left( \boldsymbol{X}
    - \E[\boldsymbol{X}] , \hat{\boldsymbol{X}}
    - \E [\boldsymbol{X}] \right)\mathrm{d} t\Big],
\end{aligned}
\end{equation}
where $\|\bm{X}\|^2\coloneqq \int_0^T |\bm{X}|_2^2\,\text{d}t$,
$|\cdot|_2$ denotes the $\ell^2$ norm of a vector, and $(\cdot,
\cdot)$ is the inner product of two $d$-dimensional vectors.  In view
of the independence between $\boldsymbol{X}- \E[\boldsymbol{X}]$ and
$\hat{\boldsymbol{X}} - \E [\boldsymbol{X}]$, we have
$\E\big[\big(\boldsymbol{X} - \E[\boldsymbol{X}], \hat{\boldsymbol{X}}
  - \E[\boldsymbol{X}] \big)\big] = \E\big[\left(\boldsymbol{X}-
  \E[\boldsymbol{X}] \right)\big] \cdot
\E\big[(\hat{\boldsymbol{X}} - \E[\boldsymbol{X}]
  )\big]=0$, and

\begin{equation}
  \E\big[ \| \boldsymbol{X} - \hat{\boldsymbol{X}}\|^2\big] \geq
  \E\big[\left\|\boldsymbol{X}- \E[\boldsymbol{X}]\right\|^2\big].
\end{equation}

\noindent Therefore, the optimal $\hat{\bm{X}}$ that minimizes the MSE is
$\hat{\bm{X}} = \E[\bm{X}]$, which indicates that the MSE tends
to fit the mean process $\E[\bm{X}]$ and make noise diminish. This is
not desirable when one wishes to fit a nonzero $\sigma$ in
Eq.~\eqref{SDE_representation}.


The KL divergence, in some cases, will diverge and thus is not
suitable for being used as a loss function. Here, we provide a simple
intuitive example when the KL divergence fail. If we consider the
degenerate case when $\text{d} X(t) = \d t, \text{d} \hat{X}(t) =
(1-\epsilon)\d t, t\in[0, T]$, then $D_{KL}(\mu, \hat{\mu})=\infty$ no
matter how small $\epsilon\neq0$ is because $\mu, \hat{\mu}$ has
different and degenerate support. However, from
Theorem~\ref{theorem1}, $\lim\limits_{\epsilon\rightarrow 0}W_2(\mu,
\hat{\mu})= 0$. Therefore, the KL divergence cannot effectively
measure the similarity between $\mu, \hat{\mu}$. Overall, the squared
$W_2$ distance is a better metric than some of the commonly used loss
metrics such as the MSE or the KL divergence.

\section{Proof to Theorem~\ref{theorem3}}
\label{proof_theorem3}
Here, we shall prove Theorem~\ref{theorem3}. We denote
\begin{equation}
  \Omega_N\coloneqq \{\bm{Y}(t)| \bm{Y}(t)
  =\bm{Y}(t_i)\,\,\, t\in [t_i, t_{i+1}), i<N-1;
    \,\,\bm{Y}(t) = \bm{Y}(t_i), \,\,\,t\in[t_i, t_{i+1}]\}
\end{equation}
to be the space of piecewise functions. We also define the space

\begin{equation}
  \tilde{\Omega}_N\coloneqq \{\bm{Y}_1(t)
  + \bm{Y}_2(t), \bm{Y}_1\in C([0, T]; \mathbb{R}^d), \bm{Y}_2\in \Omega_N\}.
\end{equation}
$\tilde{\Omega}_N$ is also a separable metric space because both
$\big(C([0, T]; \mathbb{R}^d), \|\cdot\|\big)$ and $\big(\Omega_N,
\|\cdot\|\big)$ are separable metric spaces. Furthermore, both the
embedding mapping from $C([0, T]; \mathbb{R}^d)$ to $\tilde{\Omega}_N$
and the embedding mapping from $\Omega_N$ to $\tilde{\Omega}_N$
preserves the $\|\cdot\|$ norm. Then, the two embedding mappings are
measurable, which enables us to define the measures on
$\mathcal{B}(\tilde{\Omega}_N)$ induced by the measures $\mu,
\hat{\mu}$ on $\mathcal{B}\big(C([0, T]; \mathbb{R}^d)\big)$ and the
measures $\mu_N, \hat{\mu}_N$ on $\mathcal{B}(\Omega_N)$. For
notational simplicity, we shall still denote those induced measures by
$\mu, \hat{\mu}, \mu_N, \hat{\mu}_N$.

Therefore, the inequality~Eq.~\eqref{triangular} is a direct result of
the triangular inequality for the Wasserstein distance
\citep{clement2008elementary} because $\bm{X}, \bm{X}_N, \hat{\bm{X}},
\hat{\bm{X}}_N\in \tilde{\Omega}_N$.

Next, we shall prove Eq.~\eqref{dtbound} when $\bm{X}(t),
\hat{\bm{X}}(t)$ are solutions to SDEs Eq.~\eqref{SDE_representation}
and Eq.~\eqref{approximate_sde}. Because $\bm{X}_N(t)$ is the
projection to $\bm{X}(t)$, the squared $W_2^2(\mu, \mu_N)$ can be
bounded by
\begin{equation}
\begin{aligned}
  W_2^2(\mu, \mu_N)&\leq \sum_{i=1}^N \int_{t_{i-1}}^{t_i}
  \E\big[\big|\bm{X}(t) - \bm{X}_N(t)\big|_2^2\big]\d t =
  \sum_{i=1}^N \int_{t_{i-1}}^{t_i} \sum_{\ell=1}^d
  \E\big[\big(X_{\ell}(t) - X_{N, \ell}(t)\big)^2\big]\d t
\end{aligned}
\end{equation}
For the first inequality above, we choose a specific
\textit{coupling}, i.e. the coupled distribution, $\pi$ of $\mu,
\mu_N$ that is essentially the ``original'' probability distribution.
To be more specific, for an abstract probability space $(\Omega,
\mathcal{A}, p)$ associated with ${\boldsymbol X}$, $\mu$ and
$\mu_{N}$ can be characterized by the \textit{pushforward} of $p$ via
${\boldsymbol X}$ and ${\boldsymbol X}_N$ respectively, i.e., $\mu =
        {\boldsymbol X}_* p$, defined by $\forall A \in
        \mathcal{B}\big(\tilde{\Omega}_N\big)$, elements in the Borel
        $\sigma$-algebra of $\tilde{\Omega}_N$,

\begin{equation}
\mu(A) = {\boldsymbol{X}}_*p(A) \coloneqq p\big({\boldsymbol X}^{-1}(A)\big),
\end{equation}
where $\boldsymbol X$ is interpreted as a measurable map from $\Omega$
to $\tilde{\Omega}_N$, and $\boldsymbol{X}^{-1}(A)$ is the preimage of
$A$ under $\boldsymbol{X}$.  Then, the coupling $\pi$ is defined by

\begin{equation}
\pi = ({\boldsymbol X}, {\boldsymbol X}_N)_* p,
\end{equation}
where $({\boldsymbol X}, {\boldsymbol X}_N)$ is interpreted as a
measurable map from $\Omega$ to $\tilde{\Omega}_N\times
\tilde{\Omega}_N$.  One can readily verify that the marginal
distributions of $\pi$ are $\mu$ and $\mu_{N}$ respectively. Recall
that $s$ represents the dimension of the standard Brownian motions in
the SDEs.\\

\noindent For each $\ell=1,...,d$, we have

\begin{equation}
\begin{aligned}
  &\sum_{i=1}^N \int_{t_{i-1}}^{t_i}  \!
  \E\big[\big(X_{\ell}(t)-X_{N, \ell}(t)\big)^2\big]\d t \\
   &\quad  \leq
(s+1)\Bigg[\sum_{i=1}^N\int_{t_{i-1}}^{t_i}
    \bigg(\E\Big[\Big(\medint\int_{t_i}^t f_{\ell}(\hat{X}(r), r)\d r\Big)^2\Big]
    + \E\Big[\Big(\medint\int_{t_i}^t
    \sum_{j=1}^s\sigma_{\ell, j}(\hat{X}(r), r)\d B_j(r)\Big)^2\Big]\bigg)\d t\Bigg]\\
  &\quad \leq (s+1)\sum_{i=1}^N
  \bigg((\Delta t)^2\E\Big[\medint\int_{t_{i-1}}^{t_i}f_{\ell}^2\d t\Big]
  + \Delta t\sum_{j}\E\Big[\medint\int_{t_{i-1}}^{t_i}\sigma_{\ell, j}^2\d t\Big]\bigg)
\end{aligned}
\end{equation}
The first inequality follows from the observation that
$\big(\sum_{i=1}^n a_i\big)^2 \leq n (\sum_{i=1}^n a_i^2)$ and
application of this observation to the integral representation of
${\bm X}(t)$. Summing over $\ell$, we have
\begin{equation}
  \Big(\sum_{i=1}^N \int_{t_{i-1}}^{t_i} \E\big[\big|\bm{X}(t)
    - \bm{X}_N(t)\big|_2^2\big]\d t\Big)^{1/2}\leq
  \sqrt{s+1}\,\Big(F(\Delta t)^2 + \Sigma\Delta t\Big)^{1/2}
\label{X_piecebound}
\end{equation}
Similarly, $W_2(\hat{\mu}, \hat{\mu}_N)$ can be bounded by

\begin{equation}
  W_2(\hat{\mu}, \hat{\mu}_N)\leq \sqrt{s+1}\, \sqrt{\hat{F}(\Delta t)^2
  + {\hat{\Sigma}}\Delta t}.
\label{Xhat_piecebound}
\end{equation}
Substituting Eq.~\eqref{X_piecebound} and Eq.~\eqref{Xhat_piecebound}
into Eq.~\eqref{triangular}, we have proved Eq.~\eqref{dtbound}. This
completes the proof of Theorem~\ref{theorem3}.

\section{Proof to Theorem~\ref{theorem4}}
\label{proof_theorem4}
We now give a proof to Theorem~\ref{theorem4}. First, we notice that
\begin{equation}
  \E\left[\big|\bm{X}(t) - \hat{\bm{X}}(t)\big|_2^2\right]\leq 2
  \big(FT+\hat{F}T+\Sigma+\hat{\Sigma}\big)<\infty, \,\,\, \forall t\in[0, T]
\label{M_condition}
\end{equation}
where $F, \hat{F}, \Sigma, \hat{\Sigma}$ are defined in
Eq.~\eqref{F_Sigma}. We denote

\begin{equation}
M\coloneqq \max_{t\in[0, T]} W_2\big(\mu(t), \hat{\mu}(t)\big)\leq 2
  \big(FT+\hat{F}T+\Sigma+\hat{\Sigma}\big).
\label{M_condition}
\end{equation}
By applying Theorem~\ref{theorem3} with $N=1$, the bound

\begin{equation}
\begin{aligned}
  &\hspace{-2mm}\inf_{\pi_i}\sqrt{\E_{\pi_i}\big[\big|\bm{X}(t_i)
      - \hat{\bm{X}}(t_i)\big|_2^2\big] \Delta t}
  - \sqrt{(s+1)\Delta t}\,\Big(\sqrt{F_i\Delta t + \Sigma_i}
  + \sqrt{\hat{F}_i\Delta t + \hat{\Sigma}_i} \Big)
  \leq  W_2(\boldsymbol{\mu}_i, \hat{\boldsymbol{\mu}}_i)\\
  &\quad\qquad\leq \inf_{\pi_i}\sqrt{\E_{\pi_i}\big[\big|\bm{X}(t_i) -
      \hat{\bm{X}}(t_i)\big|_2^2\big]\Delta t} + \sqrt{(s+1)\Delta t}\,
  \Big(\sqrt{F_i\Delta t + \Sigma_i} +
  \sqrt{\hat{F}_i\Delta t + \hat{\Sigma}_i} \Big).
\end{aligned}
\label{approx_i}
\end{equation}
holds true for all $i=1,2,...,N-1$. In Eq.~\eqref{approx_i},

\begin{equation}
\begin{aligned}
  &F_i\coloneqq\E\Big[\medint\int_{t_i}^{t_{i+1}} \sum_{\ell=1}^d
    f_{\ell}^2(\bm{X}(t),t)\d t\Big]<\infty,
\,\,\Sigma_i\coloneqq\E\Big[\medint\int_{t_i}^{t_{i+1}} \sum_{\ell=1}^d
\sum_{j=1}^s\sigma_{\ell, j}^2(\bm{X}(t),t)\d t\Big]<\infty,\\
&\hat{F}_i\coloneqq\E\Big[\medint\int_{t_i}^{t_{i+1}} \sum_{\ell=1}^d
\hat{f}_{\ell}^2(\hat{\bm{X}}(t),t)\d t\Big]<\infty,\,\,
\hat{\Sigma}_i\coloneqq\E\Big[\medint\int_{t_i}^{t_{i+1}}
\sum_{\ell=1}^d\sum_{j=1}^s\hat{\sigma}_{\ell, j}^2(\hat{\bm{X}}(t),t)\d t\Big]<\infty,
\end{aligned}
\end{equation}
which results from

\begin{equation}
  \sum_{i=0}^{N-1}F_i = F<\infty,\,\,
  \sum_{i=0}^{N-1}\hat{F}_i = \hat{F}<\infty, \,\,
  \sum_{i=0}^{N-1}\Sigma_i = \Sigma<\infty, \,\,
  \sum_{i=0}^{N-1}\hat{\Sigma}_i = \hat{\Sigma}<\infty,
\end{equation}
where $F, \hat{F}, \Sigma, \hat{\Sigma}$ are defined in Eq.~\eqref{F_Sigma}.
Squaring the inequality~\eqref{approx_i}, we have

\begin{equation}
\begin{aligned}
W_2^2(\boldsymbol{\mu}_i, \hat{\boldsymbol{\mu}}_i)\leq &
\inf_{\pi_i}\E_{\pi_i}\big[\big|\bm{X}(t_i)
  - \hat{\bm{X}}(t_i)\big|_2^2\big]\Delta t  \\
\: & \quad  + 2\inf_{\pi_i}\sqrt{\E_{\pi_i}\big[\big|\bm{X}(t_i)
    - \hat{\bm{X}}(t_i)\big|_2^2\big]}
\sqrt{(s+1) \Delta t}\left(\sqrt{F_i\Delta t + \Sigma_i} +
\sqrt{\hat{F}\Delta t + \hat{\Sigma}_i}\right) \\
\: & \quad + 2 (s+1)\Delta t \big(F_i \Delta t + \Sigma_i
+ \hat{F}_i \Delta t +\hat{\Sigma}_i\big),\\[5pt]
W_2^2(\boldsymbol{\mu}_i, \hat{\boldsymbol{\mu}}_i)\geq &
\inf_{\pi_i}\E_{\pi_i}\big[\big|\bm{X}(t_i)
  - \hat{\bm{X}}(t_i)\big|_2^2\big]\Delta t \\
\: & \quad - 2W_2(\boldsymbol{\mu}_i, \hat{\boldsymbol{\mu}}_i)\sqrt{(s+1)\Delta t}
\left(\sqrt{F_i\Delta t + \Sigma_i} +
\sqrt{\hat{F}\Delta t + \hat{\Sigma}_i} \right)\\
\: & \quad - 2(s+1)\Delta t \big(F_i \Delta t + \Sigma_i
+ \hat{F}_i \Delta t +\hat{\Sigma}_i\big)
\end{aligned}
\label{both_side}
\end{equation}
Specifically, from Eq.~\eqref{M_condition} and
Eq.~\eqref{approx_i},

\begin{equation}
  W_2(\boldsymbol{\mu}_i, \hat{\boldsymbol{\mu}}_i)\leq
  \sqrt{\Delta t}
  \left[M+\sqrt{s+1}\Big(\sqrt{FT+\Sigma}+\sqrt{\hat{F}T + \hat{\Sigma}}\Big)\right]
  \coloneqq \tilde{M}\sqrt{\Delta t}, \quad \tilde{M}<\infty
\end{equation}
Summing over $i=1,...,N-1$ for both inequalities in
Eq.~\eqref{both_side} and noting that $\Delta t=\frac{T}{N}$, we
conclude

\begin{equation}
\begin{aligned}
\sum_{i=1}^{N-1}W_2^2(\boldsymbol{\mu}_i, \hat{\boldsymbol{\mu}}_i)
  \leq & \sum_{i=1}^{N-1}\inf_{\pi_i}
  \E_{\pi_i}\left[\big|\bm{X}(t_i) - \hat{\bm{X}}(t_i)\big|_2^2\right]\Delta t \\
\: &\quad+ 2M\Delta t\sqrt{s+1}\sum_{i=1}^{N-1}
\left(\sqrt{F_i\Delta t + \Sigma_i} +
\sqrt{\hat{F}_i\Delta t + \hat{\Sigma}_i} \right) \\
\: & \quad +2\Delta t (s+1)
\big(F\Delta t + \Sigma + \hat{F}\Delta t + \hat{\Sigma}\big),\\
\leq & \sum_{i=1}^{N-1}\inf_{\pi_i}\E_{\pi_i}\left[\big|\bm{X}(t_i)
  - \hat{\bm{X}}(t_i)\big|_2^2\right]\Delta t \\
\: & \quad + 2 (s+1)\Delta t
  \big(F\Delta t + \Sigma + \hat{F}\Delta t + \hat{\Sigma}\big)\\
\: & \quad + M\sqrt{(s+1)\Delta t}\,
  \big((F+\hat{F}+2T)\sqrt{\Delta t} + \Sigma+\hat{\Sigma} + 2T\big)
\end{aligned}
\label{approx_i1}
\end{equation}
and

\begin{equation}
\begin{aligned}
\sum_{i=1}^{N-1}W_2^2(\boldsymbol{\mu}_i,
  \hat{\boldsymbol{\mu}}_i)\geq & \sum_{i=1}^{N-1}
  \inf_{\pi_i}\E_{\pi_i}\big[|\bm{X}(t_i) - \hat{\bm{X}}(t_i)|_2^2\big]\Delta t \\
\: & \quad -2\tilde{M}\Delta t\sqrt{s+1}\sum_{i=1}^{N-1}
\left(\sqrt{F_i\Delta t + \Sigma_i} +
\sqrt{\hat{F}_i\Delta t + \hat{\Sigma}_i} \right) \\
\: & \quad  -2 (s+1)\Delta t\,
\big(F\Delta t + \Sigma + \hat{F}\Delta t + \hat{\Sigma}\big), \\
\geq & \sum_{i=1}^{N-1}\inf_{\pi_i}
\E_{\pi_i}\big[|\bm{X}(t_i)- \hat{\bm{X}}(t_i)|_2^2\big]\Delta t \\
\: & \quad  - 2(s+1)\Delta t\,
\big(F\Delta t + \Sigma + \hat{F}\Delta t + \hat{\Sigma}\big)\\
\: & \quad - \tilde{M}\sqrt{(s+1)\Delta t}\,
  \big((F+\hat{F}+2T)\sqrt{\Delta t} + \Sigma+\hat{\Sigma} + 2T\big).
\end{aligned}
\label{approx_i2}
\end{equation}
Eqs.~\eqref{approx_i1} and \eqref{approx_i2} indicate that as
$N\rightarrow\infty$,

\begin{equation}
  \sum_{i=1}^{N-1}\inf_{\pi_i}\E_{\pi_i}\left[\big|\bm{X}(t_i) -
    \hat{\bm{X}}(t_i)\big|_2^2\right]\Delta t
  -  \sum_{i=1}^{N-1}W_2^2(\boldsymbol{\mu}_i, \hat{\boldsymbol{\mu}}_i)\rightarrow 0,
\end{equation}

Suppose $0=t_0^1<t_1^1<...<t_{N_1}^1=T$;
$0=t_0^2<t_1^2<...<t_{N_2}^2=T$ to be two sets of grids on $[0,
  T]$. We define a third set of grids $0=t_0^3<...<t_{N_3}^3=T$ such
that $\{t_0^1,...,t_{N_1}^1\}\cup \{t_0^2,...,t_{N_2}^2\} =
\{t_0^3,...,t_{N_3}^3\}$. Let $\delta
t\coloneqq\max\{\max_i(t_{i+1}^1-t_i^1), \max_j(t_{j+1}^2-t_j^2),
\max_k(t_{k+1}^3-t_k^3)\}$. We denote $\mu(t_i^1)$ and
$\hat{\mu}(t_i^1)$ to be the probability distribution of
$\bm{X}(t_i^s)$ and $\hat{\bm{X}}(t_i^s)$, $s=1,2,3$, respectively.
We will prove

\begin{equation}
  \Big|\sum_{i=0}^{N_1-1}W_2^2\big(\mu(t_i^1), \hat{\mu}(t_i^1)\big)(t_{i+1}^1-t_i^1)
  - \sum_{i=0}^{N_3-1}W_2^2\big(\mu(t_i^3),
  \hat{\mu}(t_i^3)\big)(t_{i+1}^3-t_i^3)\Big| \rightarrow 0,
\label{limit_exist}
\end{equation}
as $\delta t\rightarrow 0$.

First, suppose in the interval $(t_i^1, t_{i+1}^1)$, we have
$t_i^1=t_{\ell}^3<t_{\ell+1}<...<t_{\ell+s}^3=t_{i+1}^1, s\geq 1$,
then for $s>1$, since
$t_{i+1}^1-t_i^1=\sum_{k=\ell}^{\ell+s-1}(t_{k+1}^3-t_k^3)$, we have

\begin{equation}
\begin{aligned}
  &\Big|W_2^2\big(\mu(t_i^1), \hat{\mu}(t_i^1)\big)(t_{i+1}^1-t_i^1)
  - \sum_{k=\ell}^{\ell+s-1}
  W_2^2\big(\mu(t_k^3)), \hat{\mu}(t_i^3)\big)(t_{k+1}^3-t_k^3)\Big| \qquad \:\\
  &\quad\quad\leq\sum_{k=\ell+1}^{\ell+s-1}
  \Big|W_2\big(\mu(t_i^1), \hat{\mu}(t_i^1)\big)
  - W_2\big(\hat{\mu}(t_i^3), \hat{\mu}(t_k^3)\big)\Big|\\
  &\quad\quad\qquad\qquad\times\Big(W_2\big(\mu(t_i^1),
  \hat{\mu}(t_i^1)\big)+ W_2\big(\mu(t_k^3),
  \hat{\mu}(t_k^3)\big)\Big) (t_{k+1}^3 - t_k^3).
\end{aligned}
\label{triang}
\end{equation}
On the other hand, because we can take a specific coupling
  $\pi^*$ to be the joint distribution of $(\bm{X}(t_i^1),
  \bm{X}(t_k^3))$,
\begin{equation}
\begin{aligned}
  &W_2(\mu(t_i^1), \mu(t_k^3)) \leq
  \Big(\E\big[|\bm{X}(t_k^3) - \bm{X}(t_i^1)|_2^2\big]\Big)^{1/2}\\
  &\quad\quad\leq \sqrt{s+1}\,\, \E\Big[\int_{t_i}^{t_{i+1}} \sum_{\ell=1}^d
    f_{\ell}^2(\bm{X}(t),t) \d t+ \sum_{\ell=1}^d
\sum_{j=1}^s\sigma_{\ell, j}^2(\bm{X}(t),t)\d t\Big]^{1/2}.
\end{aligned}
\end{equation}
Similarly, we have

\begin{equation}
  W_2\big(\hat{\mu}(t_i^1), \hat{\mu}(t_k^3)\big)
  \leq \sqrt{s+1}\,\, \E\Big[\int_{t_i}^{t_{i+1}} \sum_{\ell=1}^d
    \hat{f}_{\ell}^2(\bm{X}(t),t) \d t + \sum_{\ell=1}^d
\sum_{j=1}^s\hat{\sigma}_{\ell, j}^2(\bm{X}(t),t)\d t\Big]^{1/2}
\end{equation}
From the triangular inequality of the Wasserstein distance, we find
\begin{equation}
  \Big|W_2\big(\mu(t_i^1), \hat{\mu}(t_i^1)\big)-
  W_2\big(\mu(t_k^3), \hat{\mu}(t_k^3)\big)\Big|
  \leq W_2\big(\mu(t_i^1), \mu(t_k^3)\big)
  + W_2\big(\hat{\mu}(t_i^1), \hat{\mu}(t_k^3)\big).
\label{traing_property}
\end{equation}
Substituting Eq.~\eqref{traing_property} into Eq.~\eqref{triang}, we
conclude that

\begin{equation}
\begin{aligned}
  &\Big|W_2^2\big(\mu(t_i^1), \hat{\mu}(t_i^1)\big)(t_{i+1}^1-t_i^1)
  - \sum_{k=\ell}^{\ell+s-1} W_2^2\big((\mu(t_k^3),
  \hat{\mu}(t_k^3)\big)(t_{k+1}^3-t_k^3)\Big|\\
  &\qquad\qquad\qquad\quad \leq 2M(t_{i+1}^1-t_i^1)\left(\sqrt{F_i\delta t+\Sigma_i}
  + \sqrt{\hat{F}_i \delta t+ \hat{\Sigma}_i}\right).
\end{aligned}
\label{intermediate}
\end{equation}

When the conditions in Eq.~\eqref{condition_dt} hold true, we use
Eq.~\eqref{intermediate} in Eq.~\eqref{limit_exist} to find
\begin{equation}
\begin{aligned}
  &\Big|\sum_{i=0}^{N_1-1}W_2^2\big(\mu(t_i^1), \hat{\mu}(t_i^1)\big)(t_{i+1}^1-t_i^1)
  - \sum_{i=0}^{N_3-1}W_2^2\big(\mu(t_i^3), \hat{\mu}(t_i^3)\big)(t_{i+1}^3-t_i^3)\Big| \\
  &\qquad\qquad\qquad\quad \leq 2M T\max_i\left(\sqrt{F_i\delta t+\Sigma_i}
  + \sqrt{\hat{F}_i\delta t + \hat{\Sigma}_i}\right)\rightarrow0
\end{aligned}
\end{equation}
as $\delta t\rightarrow 0$. Similarly,

\begin{equation}
\begin{aligned}
  &\Big|\sum_{i=0}^{N_2-1}W_2^2\big(\mu(t_i^2), \hat{\mu}(t_i^2)\big)(t_{i+1}^2-t_i^2)
  - \sum_{i=0}^{N_3-1}W_2^2\big(\mu(t_i^3), \hat{\mu}(t_i^3)\big)(t_{i+1}^3-t_i^3)\Big| \\
  &\qquad\qquad\qquad\quad \leq 2M T\max_i\left(\sqrt{F_i\delta t+\Sigma_i}
  + \sqrt{\hat{F}_i\delta t + \hat{\Sigma}_i}\right)\rightarrow 0
\end{aligned}
\end{equation}
as $\delta t\rightarrow 0$. Thus,

\begin{equation}
\begin{aligned}
  &\Big|\sum_{i=0}^{N_1-1}W_2^2\big(\mu(t_i^1), \hat{\mu}(t_i^1)\big)(t_{i+1}^1-t_i^1)
  - \sum_{i=0}^{N_2-1}W_2^2\big(\mu(t_i^2),\hat{\mu}(t_i^2)\big)
  (t_{i+1}^2-t_i^2)\Big|\rightarrow 0
\end{aligned}
\end{equation}
as $\delta t\rightarrow0$, which implies the limit

\begin{equation}
  \lim\limits_{N\rightarrow \infty} \sum_{i=1}^{N-1}\inf_{\pi_i}
  \E_{\pi_i}\left[\big|\bm{X}(t_i^1) -
    \hat{\bm{X}}(t_i^1)\big|_2^2\right](t_i^1-t_{i-1}^1)
  =\lim\limits_{N\rightarrow \infty}\sum_{i=1}^{N-1}
  W_2^2\big(\mu(t_i^1), \hat{\mu}(t_i^1)\big)(t_i^1-t_{i-1}^1)
\end{equation}
exists. This completes the proof of Theorem~\ref{theorem4}.

\section{Definition of different loss metrics used in the examples}
\label{def_loss}
Six loss functions for 1D cases were considered:
\begin{compactenum}
\item The squared Wasserstein-2 distance~(Eq.~\eqref{time_discretize})
$$W_2^2(\mu_N^{\text{e}}, \hat{\mu}_N^{\text{e}}),$$ where
$\mu_N^{\text{e}}$ and $\hat{\mu}_N^{\text{e}}$ are the empirical
distributions of the vector $({X}(t_1),...X(t_{N-1}))$ and
$(\hat{{X}}(t_1),...,\hat{{X}}(t_{N-1}))$, respectively. It is
estimated by
\begin{equation}
W_2^2(\mu_N^{\text{e}},
\hat{\mu}_N^{\text{e}})\approx\texttt{ot.emd2}\Big(\frac{1}{M}\bm{I}_{M},
\frac{1}{M}\bm{I}_{M}, \bm{C}\Big),
\label{time_coupling}
\end{equation}
where $\texttt{ot.emd2}$ is the function for solving the earth movers
distance problem in the $\texttt{ot}$ package of Python, $M$ is the
number of ground-truth and predicted trajectories, $\bm{I}_{\ell}$ is
an $M$-dimensional vector whose elements are all 1, and
$\bm{C}\in\mathbb{R}^{M\times M}$ is a matrix with entries
$(\bm{C})_{ij} = (X_N^i-\hat{X}^j_N)_2^2$. $X^i _N$ is the
vector of the values of the $i^{\text{th}}$ ground-truth trajectory at
time points $t_1,...,t_{N-1}$, and $\hat{X}^j _N$ is the vector
of the values of the $j^{\text{th}}$ predicted trajectory at time
points $t_1,...,t_{N-1}$.
\item The squared time-decoupled Wasserstein-2 distance averaged over
  each time step (Eq.~\eqref{approximation}):
  $$\tilde{W}^2_2(\mu_N, \hat{\mu_N}) = \sum_{i=1}^{N-1}
  W^2_2(\mu_{N}^{\text{e}}(t_i),
  \hat{\mu}^{\text{e}}_{N}(t_i))\Delta{t}$$, where $\Delta t$ is the
  time step and $W_2$ is the Wasserstein-2 distance between two
  empirical distributions $\mu_{N}^{\text{e}}(t_i),
  \hat{\mu}_{N}^{\text{e}}(t_i)$.  These distributions are calculated
  by the samples of the trajectories of $X(t), \hat{X}(t)$ at a given
  time step $t=t_i$, respectively.
\item Mean squared error (MSE) between the trajectories, where $M$ is
  the total number of the ground-truth and prediction
  trajectories. $X_{i, j}$ and $\hat{X}_{i, j}$ are the values of the
  $j^{\text{th}}$ ground-truth and prediction trajectories at time
  $t_i$, respectively:
$$\operatorname{MSE}({X}, \widehat { X}) = \sum_{i=1}^N \sum_{j=1}^M
(X_{i, j}-\hat X_{i, j})^2\Delta{t}.$$
\item The sum of squared distance between mean trajectories and
  absolute distance between trajectories, which is a common practice
  for estimating the parameters of an SDE. Here $M$ and $X_{i, j}$ and
  $\hat{X}_{i, j}$ have the same meaning as in the MSE
  definition. $\operatorname{var}(X_i)$ and $\operatorname{var}
  (\hat{X}_i)$ are the variances of the empirical distributions of
  $X(t_i), \hat{X}(t_i)$, respectively. We shall denote this loss
  function by
$$(\operatorname{mean}^2+\operatorname{var})({ X}, \hat { X}) =
  \sum_{i=1}^N \left[\Big( \frac{1}{n}\sum_{j=1}^M X_{i, j} -
    \frac{1}{n}\sum_{i=1}^N \hat{X}_{i, j} \Big)^2 + \left|
    \operatorname{var}(X_i) - \operatorname{var} (\hat{X}_i)
    \right|\right]\Delta{t}.$$
\item Negative approximate log-likelihood of the trajectories:
$$ - \log \mathcal{L}( X| \sigma) = - \sum_{i=0}^{N-1} \sum_{j=1}^{M}
  \log \rho_\mathcal{N}\Big[\frac{X_{i+1,j} -{X}_{t,j}+ f(X_{i, j},
      t_i) \Delta t }{ \sigma^2(X_{i, j}, t_i) \Delta t}\Big],
$$
where $\rho_{\mathcal{N}}$ is the probability density function of the
standard normal distribution and $f(X_{i, j}, t_i), \sigma(X_{i, j},
t_i)$ are the ground-truth drift and diffusion functions in
Eq.~\eqref{SDE_representation}. $M$ and $X_{i, j}$ and $\hat{X}_{i,
  j}$ have the same meaning as in the MSE definition.
\item MMD (maximum mean discrepancy) 
  \citep{li2015generative}:
$$
\text{MMD}(X, \hat{X}) = \sum_{i=1}^N\big(\E_p[K(X_i, X_i)]
- 2\E_{p, q}[K(X_i, \hat{X}_i)] + \E_q[K(\hat{X}_i, \hat{X}_i)]\big)\Delta{t},
$$
where $K$ is the standard radial basis function (or Gaussian kernel)
with multiplier $2$ and number of kernels $5$. $X_i$ and $\hat{X}_i$
are the values of the ground-truth and prediction trajectories at time
$t_j$, respectively.

\end{compactenum}

\noindent Five $W_2$ distance based loss functions for the 2D SDE reconstruction
problem Example~\ref{example4} are listed as follows

\begin{compactenum}
\item 2D squared $W_2$ loss
$$ \sum_{i=1}^{N-1} \Big(W_2^2\big(\mu_{N, 1}(t_i), \hat{\mu}_{N, 1}(t_i)\big) +
  W_2^2\big(\mu_{N, 2}(t_i), \hat{\mu}_{N, 2}(t_i)\big)\Big)\Delta{t}
$$ where $\mu_{N, 1}(t_i)$ and $\hat{\mu}_{N, 1}(t_i)$ are the empirical
distributions of $X_1, \hat{X}_1$ at time $t_i$, respectively. Also,
$\mu_{N, 2}(t_i)$ and $\hat{\mu}_{N, 2}(t_i)$ are the empirical
distributions of $X_2, \hat{X}_2$ at time $t_i$, respectively.
\item Weighted sliced squared $W_2$ loss
$$
\sum_{i=1}^{N-1} \Big(\sum_{k=1}^m \frac{N_k}{\sum_{\ell=1}^m N_\ell}
W_2^2\big(\mu_{N, k}^{\text{s}}(t_i), \hat{\mu}_{N, k}^{\text{s}}(t_i)\big)\Big)\Delta{t}
$$ where $\mu_{N, k}^{\text{s}}(t_i)$ is the empirical distribution for
$\sqrt{X_1(t_i)^2+X_2(t_i)^2}$ such that the angle between the two
vectors $\big(X_1(t_i), X_2(t_i)\big)$ and $(1, 0)$ is in
$[\frac{2(k-1)\pi}{m}, \frac{2k\pi}{m})$; $\hat{\mu}_{N,
    k}^{\text{s}}(t_i)$ is the empirical distribution for
  $\sqrt{\hat{X}_1(t_i)^2+\hat{X}_2(t_i)^2}$ such that the angle
  between the two vectors $\big(\hat{X}_1(t_i), \hat{X}_2(t_i)\big)$
  and $(1, 0)$ is in $[\frac{2(k-1)\pi}{m}, \frac{2k\pi}{m})$; $N_k$
    is the number of predictions such that the angle between the two
    vectors $(\hat{X}_1(t_i), \hat{X}_2(t_i))$ and $(1, 0)$ is in
    $\big[\frac{2(k-1)\pi}{m}, \frac{2k\pi}{m}\big)$.
\item The loss function~Eq.~\eqref{time_discretize}
$$W_2^2(\mu_N^{\text{e}}, \hat{\mu}_N^{\text{e}}),$$ where
$\mu_N^{\text{e}}$ and $\hat{\mu}_N^{\text{e}}$ are the empirical
distributions of the vector $(\bm{X}(t_1),...\bm{X}(t_{N-1}))$ and
$(\hat{\bm{X}}(t_1),...,\hat{\bm{X}}(t_{N-1}))$, respectively. It is
estimated by
\begin{equation}
W_2^2(\mu_N^{\text{e}},
\hat{\mu}_N^{\text{e}})\approx\texttt{ot.emd2}(\frac{1}{M}\bm{I}_{M},
\frac{1}{M}\bm{I}_{M}, \bm{C}),
\label{time_coupling}
\end{equation}
where $\texttt{ot.emd2}$ is the function for solving the earth movers
distance problem in the $\texttt{ot}$ package of Python, $M$ is the
number of ground-truth and predicted trajectories, $\bm{I}_{\ell}$ is
an $M$-dimensional vector whose elements are all 1, and
$\bm{C}\in\mathbb{R}^{M\times M}$ is a matrix with entries
$(\bm{C})_{ij} = |\bm{X}_N^i-\hat{\bm{X}}^j_N|_2^2$. $\bm{X}^i _N$ is
the vector of the values of the $i^{\text{th}}$ ground-truth
trajectory at time points $t_1,...,t_{N-1}$, and $\hat{\bm{X}}^j _N$
is the vector of the values of the $j^{\text{th}}$ predicted
trajectory at time points $t_1,...,t_{N-1}$.
\item The right-hand side of~Eq.~\eqref{approximation}. It is
estimated by
\begin{equation}
\begin{aligned}
&\sum_{i=1}^{N-1}
\inf_{\pi_i} \E_{\pi_i} [|\bm{X}(t_i) -
\hat{\bm{X}}(t_i)|^2_2] \Delta t\\[-3pt]
&\qquad\qquad\approx\sum_{i=1}^{N-1}
W_2^2\big(\mu_{N}^{\text{e}}(t_i),
\hat{\mu}_{N}^{\text{e}}(t_i)\big)\Delta t
\approx\Delta{t}\sum_{i=1}^{N-1}\texttt{ot.emd2}\Big(\frac{1}{M}\bm{I}_{M},
\frac{1}{M}\bm{I}_{M}, \bm{C}_i\Big),
\end{aligned}
\label{eq17_def}
\end{equation}
where $\mu_{N}^{\text{e}}(t_i), \hat{\mu}_{N}^{\text{e}}(t_i)$ are the
empirical distribution of $\bm{X}(t_i)$, $\hat{\bm{X}}(t_i)$,
respectively, and $\texttt{ot.emd2}$ is the function for solving the
earth movers distance problem in the $\texttt{ot}$ package of
Python. $M$ is the number of ground-truth and predicted trajectories,
and $\bm{I}_{M}$ is an $\ell$-dimensional vector whose elements are
all 1. Here, the matrix $\bm{C}_i\in\mathbb{R}^{M\times{M}}$ has
entries $(\bm{C}_i)_{sj} = |\bm{X}^s(t_i)-\hat{\bm{X}}^j(t_i)|_2^2$
for $i=1,...,N-1$. $\bm{X}^s(t_i)$ is the vector of the values of the
$s^{\text{th}}$ ground-truth trajectory at the time point $t_i$, and
$\hat{\bm{X}}^j(t_i)$ is the vector of the values of the
$j^{\text{th}}$ predicted trajectory at the time point $t_i$.
\item MMD (maximum mean discrepancy) 
\citep{li2015generative}: 
$$
\text{MMD}(\bm{X}, \hat{\bm{X}}) = \sum_{i=1}^N
\big(\E_p[K(\bm{X}_i, \bm{X}_i)]
- 2\E_{p, q}[K(\bm{X}_i, \hat{\bm{X}}_i)]
+ \E_q[K(\hat{\bm{X}}_i, \hat{\bm{X}}_i)]\big)\Delta{t},
$$
where $K$ is the standard radial basis function (or Gaussian kernel)
with multiplier $2$ and number of kernels $5$. $\bm{X}_i$ and
$\hat{\bm{X}}_i$ are the values of the ground-truth and prediction
trajectories at time $t_j$, respectively.
\end{compactenum}

\section{Default training setting}
\label{training_details}
Here we list the default training hyperparameters and gradient descent
methods for each example in Table~\ref{tab:setting}.
\begin{table}[thbp]
\centering
\begin{tabular}{lllll}
\toprule
{Loss} & Example 1 & Example 2 & Example 3 & Example 4 \\
\midrule
Gradient descent method & AdamW & AdamW & AdamW & AdamW \\
Learning rate & 0.001 & 0.002 & 0.002 & 0.0005 \\
Weight decay & 0.005 & 0.005 & 0.005 & 0.005 \\
Number of epochs & 1000 & 2000 & 2000 & 2000 \\
Number of samples & 100 & 200 & 256 & 200 \\
Hidden layers in $\Theta_1$ & 2 & 1 & 1 & 1 \\
Neurons in each layer in $\Theta_1$ & 32 & 32 & 32 & 32 \\
Hidden layers in $\Theta_2$ & 2 & 1 & 1 & 1 \\
Activation function & tanh &ReLu &ReLu &ReLu \\
Neurons in each layer in $\Theta_2$ & 32 & 32 & 32 & 32 \\
$\Delta t$ & 0.1 & 0.05 & 1 & 0.02 \\
\bottomrule
\end{tabular}
\caption{Training settings for each example.} 
\label{tab:setting}
\end{table}

\section{Uncertainty in the initial condition}
\label{IC_sensitivity}
For reconstructing the CIR model Eq.~\eqref{CIRmodel} in
Example~\ref{example2}, instead of using the same initial condition
for all trajectories, we shall investigate the numerical performance
of our proposed squared $W_2$ distance loss when the initial condition
is not fixed, but rather sampled from a distribution.

First, we construct an additional dataset of the CIR model to allow
the initial value $u_0 \sim \mathcal{N}(2, \delta^2)$, with $\delta^2$
ranging from 0 to 1, and $\mathcal{N}$ stands for the 1D normal
distribution. We then train the model by minimizing
Eq.~\eqref{approximation} to reconstruct Eq.~\eqref{CIRmodel} with the
same hyperparameters as in Example~\ref{example2}. The results are
shown in Table~\ref{tab:cir_delta}, which indicate our proposed
squared $W_2$ loss function is rather insensitive to the ``noise",
\textit{i.e.}, the variance in the distribution of the initial
condition.

\begin{table}[thbp]
\centering
\begin{tabular}{lcllr}
\toprule
{Loss} & \(\delta\) & {{Relative Errors in \(f\)}} & {{Relative Errors in \(\sigma\)}} & {\(N_{\rm repeats}\)} \\
\midrule
$W_2$ & 0.0 & 0.072 ($\pm$ 0.008) & 0.071 ($\pm$ 0.023) & 10 \\
$W_2$ & 0.1 & 0.053 ($\pm$ 0.008) & 0.043 ($\pm$ 0.016) & 10 \\
$W_2$ & 0.2 & 0.099 ($\pm$ 0.007) & 0.056 ($\pm$ 0.019) & 10 \\
$W_2$ & 0.3 & 0.070 ($\pm$ 0.014) & 0.083 ($\pm$ 0.026) & 10 \\
$W_2$ & 0.4 & 0.070 ($\pm$ 0.014) & 0.078 ($\pm$ 0.040) & 10 \\
$W_2$ & 0.5 & 0.075 ($\pm$ 0.013) & 0.138 ($\pm$ 0.021) & 10 \\
$W_2$ & 0.6 & 0.037 ($\pm$ 0.018) & 0.069 ($\pm$ 0.017) & 10 \\
$W_2$ & 0.7 & 0.075 ($\pm$ 0.016) & 0.043 ($\pm$ 0.014) & 10 \\
$W_2$ & 0.8 & 0.041 ($\pm$ 0.012) & 0.079 ($\pm$ 0.023) & 10 \\
$W_2$ & 0.9 & 0.082 ($\pm$ 0.015) & 0.108 ($\pm$ 0.033) & 10 \\
$W_2$ & 1.0 & 0.058 ($\pm$ 0.024) & 0.049 ($\pm$ 0.025) & 10 \\
\bottomrule
\end{tabular}
\caption{Reconstructing the CIR model Eq.~\eqref{CIRmodel} when
  $u_0\sim\mathcal{N}(2, \delta^2)$ with different variance
  $\delta^2$. The results indicate that the reconstruction results are
  not sensitive to the variance in the distribution of the initial
  value $u_0$.}
\label{tab:cir_delta}
\end{table}

\section{Neural network structure}
\label{neural_structure}
We examine how the neural network structure affects the reconstruction
of the CIR model Eq.~\eqref{CIRmodel} in Example~\ref{example2}. We vary
the number of layers and the number of neurons in each layer (the
number of neurons are set to be the same in each hidden layer), and
the results are shown in Table~\ref{tab:cir_width}.

\begin{table}\centering
  \caption{Reconstructing the CIR model when using neuron networks of
    different widths and numbers in each hidden layer to parameterize
    $\hat{f}, \hat{\sigma}$ in Eq.~\eqref{approximate_sde}.}
\begin{tabular}{lrcllr}
\toprule Loss & Width & Layer & Relative Errors in $f$ & Relative
Errors in $\sigma$ & $N_{\rm repeats}$ \\ \midrule $W_2$ & 16 & 1 & \(
0.131 (\pm 0.135) \) & \( 0.170 (\pm 0.102) \) & 10 \\ $W_2$ & 32 & 1
& \( 0.041 (\pm 0.008) \) & \( 0.109 (\pm 0.026) \) & 10 \\ $W_2$ & 64
& 1 & \( 0.040 (\pm 0.008) \) & \( 0.104 (\pm 0.019) \) & 10 \\ $W_2$
& 128 & 1 & \( 0.040 (\pm 0.008) \) & \( 0.118 (\pm 0.019) \) & 10
\\ $W_2$ & 32 & 2 & \( 0.049 (\pm 0.015) \) & \( 0.123 (\pm 0.020) \)
& 10 \\ $W_2$ & 32 & 3 & \( 0.094 (\pm 0.013) \) & \( 0.166 (\pm
0.041) \) & 10 \\ $W_2$ & 32 & 4 & \( 0.124 (\pm 0.020) \) & \( 0.185
(\pm 0.035) \) & 10 \\ $W_2$ & 32 & 5 & \( 0.041 (\pm 0.008) \) & \(
0.122 (\pm 0.024) \) & 10 \\ $W_2$ & 32 & 6 & \( 0.043 (\pm 0.013) \)
& \( 0.117 (\pm 0.024) \) & 10 \\ $W_2$ & 32 & 7 & \( 0.044 (\pm
0.012) \) & \( 0.109 (\pm 0.017) \) & 10 \\ \bottomrule
\end{tabular}
\label{tab:cir_width}
\end{table}

The results in Table~\ref{tab:cir_width} show that increasing the
number of neurons in each layer improves the reconstruction accuracy
in $\sigma$. For the reconstructing CIR model in
Example~\ref{example2}, using 32 neurons in each layer seems to be
sufficient. On the other hand, when each layer contains 32 neurons,
the number of hidden layers in the neural network seems does not
affect the reconstruction accuracy of $f, \sigma$, and this indicates
even 1 or 2 hidden layers are sufficient for the reconstruction of $f,
\sigma$. Thus, reconstructing the CIR model in Example~\ref{example2}
using our proposed squared $W_2$ based loss function does not require
using complex deep or wide neural networks.

We also consider using the ResNet neural network structure
\citep{he2016deep}. However, the application of the ResNet technique
does not improve the reconstruction accuracy of the CIR model in
Example~\ref{example2}. This is because simple feedforward multilayer
neural network structure can work well for learning
Eq.~\eqref{CIRmodel} when reconstructing both $f$ and $\sigma$ so we do
not need deep neural networks. Thus. the ResNet technique is not
required. The results are shown in Table~\ref{tab:cir_resnet_layer}.

\begin{table}
  \centering
  \caption{Reconstructing the CIR model~Eq.~\eqref{CIRmodel} when neuron
    networks have different numbers of hidden layers and are equipped
    with the ResNet technique. Each hidden layer contains 32 neurons.}

\begin{tabular}{lcccc}
\toprule
Loss & Layer & Relative Errors in \( f \) & Relative Errors in \( \sigma \) & $N_{\rm repeats}$ \\
\midrule
$W_2$ & 1 & $0.045 (\pm 0.012)$ & $0.116 (\pm 0.025)$ & 10 \\
$W_2$ & 2 & $0.053 (\pm 0.011)$ & $0.108 (\pm 0.024$ & 10 \\
$W_2$ & 3 & $0.071 (\pm 0.017)$ & $0.117 (\pm 0.040)$ & 10 \\
$W_2$ & 4 & $0.096 (\pm 0.035)$ & $0.149 (\pm 0.064)$ & 10 \\
\bottomrule
\end{tabular}
\label{tab:cir_resnet_layer}
\end{table}

\section{Using the stochastic gradient descent method for optimization}
\label{sgdtest}

Here, we shall reconstruct the OU process Eq.~\eqref{OUprocess} in
Example~\ref{example3} with the initial condition $X(0)=0$ using the
MMD and our squared $W_2$ distance loss functions
Eqs.~\eqref{time_discretize} and~\eqref{approximation} with different
numbers of ground-truth trajectories and different batch sizes for
applying the stochastic gradient descent technique for optimizing the
parameters in the neural networks for reconstructing the SDE.
{\tiny\begin{table}[htbp]
\caption{Errors and runtime for different loss functions and different
  numbers of ground-truth trajectories when the training batch size is
  fixed to 16 and 256. The MMD and our proposed squared $W_2$ distance
  Eq.~\eqref{time_discretize} and well as our proposed time-decoupled
  squared $W_2$ distance Eq.~\eqref{approximation} are used as the
  loss function.}
\label{tab:performance_metrics}
\centering
\begin{tabular}{ccccccccccc}
\hline
Loss & $N_{\text{sample}}$ & Batch size & Relative error in $f$  & Relative error in $\sigma$ & Runtime  & $N_{\text{repeats}}$ \\
\hline
MMD & 64 & 16 & 0.30 $\pm$ 0.12 & 0.49 $\pm$ 0.17 & 1.19 $\pm$ 0.59 & 10 \\
MMD & 128 & 16 & 0.30 $\pm$ 0.09 & 0.50 $\pm$ 0.20 & 1.27 $\pm$ 0.58 & 10 \\
MMD & 256 & 16 & 0.31 $\pm$ 0.09 & 0.44 $\pm$ 0.21 & 1.31 $\pm$ 0.59 & 10 \\
MMD & 512 & 16 & 0.22 $\pm$ 0.12 & 0.43 $\pm$ 0.18 & 1.22 $\pm$ 0.37 & 10 \\
MMD & 1024 & 16 & 0.23 $\pm$ 0.11 & 0.37 $\pm$ 0.24 & 1.70 $\pm$ 0.47 & 10 \\
Eq.~\eqref{time_discretize} & 64 & 16 & 0.28 $\pm$ 0.06 & 0.66 $\pm$ 0.11 & 0.83 $\pm$ 0.26 & 10 \\
Eq.~\eqref{time_discretize} & 128 & 16 & 0.24 $\pm$ 0.07 & 0.68 $\pm$0.11 & 0.73 $\pm$ 0.18 & 10 \\
Eq.~\eqref{time_discretize} & 256 & 16 & 0.25 $\pm$ 0.07 & 0.66 $\pm$ 0.09 & 0.67 $\pm$ 0.14 & 10 \\
Eq.~\eqref{time_discretize} & 512 & 16 & 0.23 $\pm$ 0.06 & 0.68 $\pm$ 0.09 & 0.75 $\pm$ 0.16 & 10 \\
Eq.~\eqref{time_discretize} & 1024 & 16 & 0.25 $\pm$ 0.07 & 0.66 $\pm$ 0.09 & 1.02 $\pm$ 0.47 & 10 \\
Eq.~\eqref{approximation}  & 64 & 16 & 0.20 $\pm$ 0.06 & 0.42 $\pm$ 0.08 & 0.61 $\pm$ 0.14 & 10 \\
Eq.~\eqref{approximation} & 128 & 16 & 0.22 $\pm$ 0.06 & 0.37 $\pm$ 0.14 & 0.78 $\pm$ 0.35 & 10 \\
Eq.~\eqref{approximation}  & 256 & 16 & 0.21 $\pm$ 0.07 & 0.39 $\pm$ 0.16 & 0.88 $\pm$ 0.46 & 10 \\
Eq.~\eqref{approximation} & 512 & 16 & 0.23 $\pm$ 0.06 & 0.43 $\pm$ 0.15 & 0.72 $\pm$ 0.11 & 10 \\
Eq.~\eqref{approximation} & 1024 & 16 & 0.21 $\pm$ 0.03 & 0.36 $\pm$ 0.12 & 1.08 $\pm$ 0.52 & 10 \\
MMD & 64 & 256 & 0.26 $\pm$ 0.12 & 0.41 $\pm$ 0.20 & 1.54 $\pm$ 0.66 & 10 \\
MMD & 128 & 256 & 0.25 $\pm$ 0.14 & 0.40 $\pm$ 0.23 & 1.82 $\pm$ 0.78 & 10 \\
MMD & 256 & 256 & 0.25 $\pm$ 0.12 & 0.35 $\pm$ 0.21 & 3.68 $\pm$ 1.31 & 10 \\
MMD & 512 & 256 & 0.23 $\pm$ 0.14 & 0.37 $\pm$ 0.23 & 3.45 $\pm$ 1.50 & 10 \\
MMD & 1024 & 256 & 0.23 $\pm$ 0.13 & 0.35 $\pm$ 0.21 & 3.09 $\pm$ 1.35 & 10 \\
Eq.~\eqref{time_discretize} & 64 & 256 & 0.28 $\pm$ 0.08 & 0.61 $\pm$ 0.04 & 1.19 $\pm$ 0.45 & 10 \\
Eq.~\eqref{time_discretize} & 128 & 256 & 0.31 $\pm$ 0.07 & 0.61 $\pm$ 0.07 & 1.04 $\pm$ 0.48 & 10 \\
Eq.~\eqref{time_discretize} & 256 & 256 & 0.26 $\pm$ 0.07 & 0.53 $\pm$ 0.03 & 0.96 $\pm$ 0.43 & 10 \\
Eq.~\eqref{time_discretize} & 512 & 256 & 0.26 $\pm$ 0.08 & 0.56 $\pm$ 0.05 & 0.98 $\pm$ 0.40 & 10 \\
Eq.~\eqref{time_discretize} & 1024 & 256 & 0.27 $\pm$ 0.08 & 0.56 $\pm$ 0.05 & 0.89 $\pm$ 0.36 & 10 \\
Eq.~\eqref{approximation} & 64 & 256 & 0.24 $\pm$ 0.08 & 0.41 $\pm$ 0.13 & 1.39 $\pm$ 0.53 & 10 \\
Eq.~\eqref{approximation}  & 128 & 256 & 0.26 $\pm$ 0.11 & 0.37 $\pm$ 0.17 & 1.36 $\pm$ 0.61 & 10 \\
Eq.~\eqref{approximation}  & 256 & 256 & 0.20 $\pm$ 0.08 & 0.31 $\pm$ 0.16 & 1.72 $\pm$ 0.73 & 10 \\
Eq.~\eqref{approximation}  & 512 & 256 & 0.25 $\pm$ 0.11 & 0.38 $\pm$ 0.20 & 1.67 $\pm$ 0.73 & 10 \\
Eq.~\eqref{approximation}  & 1024 & 256 & 0.26 $\pm$ 0.10 & 0.39 $\pm$ 0.20 & 1.64 $\pm$ 0.79 & 10 \\
\hline
\end{tabular}
\end{table}}

We train 2000 epochs with a learning rate 0.001 for all numerical
experiments. In all cases, the loss functions converge before 2000
epochs. From Table~\ref{tab:performance_metrics}, for all three loss
function, \textit{i.e.}, the MMD loss, Eq.~\eqref{time_discretize},
and Eq.~\eqref{approximation}, a larger number of training samples
leads to more accurate reconstruction of $\sigma$ (the noise
term). Furthermore, it can be seen from
Table~\ref{tab:performance_metrics} that using a smaller batch size
(16) for training tends to lead to less accurate reconstruction of
$\sigma$ for the MMD and Eq.~\eqref{time_discretize} loss functions
even if the number of trajectories in the training set is large. This
feature might arise because the trajectories are intrinsically noisy
and evaluating MMD and Eq.~\eqref{time_discretize} will be inaccurate
if the batch size is small.  Therefore, using a smaller batch size
does not remedy the high cost of MMD as the reconstruction error is
large and leads to inaccurate reconstruction of the ground-truth SDE
for smaller $N_{\text{sample}}$.  On the other hand, our proposed
time-decoupled squared $W_2$ distance loss function
Eq.~\eqref{approximation} gives similar performance in reconstructing
$f, \sigma$ for both a batch size of $16$ and a batch size of 256. In
other words, using Eq.~\eqref{approximation} is more robust to a
smaller batch size.  From Table~\ref{tab:performance_metrics}, using a
smaller batch size (16) leads to faster training. Thus, we can
consider using Eq.~\eqref{approximation} as the loss function together
with a smaller batch size to boost training efficiency.

From the results in both Example~\ref{example3} and
Table~\ref{tab:performance_metrics}, our proposed time-decoupled
squared $W_2$ distance Eq.~\eqref{approximation} is faster and more
efficient than the MMD method and Eq.~\eqref{time_discretize}, making it
potentially most suitable among all three loss functions for
reconstructing SDEs.

\section{Additional discussion on the loss functions
  Eqs.~\eqref{time_discretize} and~\eqref{approximation}}
\label{appendix_more_discussion}
Here, we make an additional comparison between using
Eq.~\eqref{time_discretize} and Eq.~\eqref{approximation} as loss
functions in Example 4. We set the number of training samples to be
128 and other hyperparameters for training to be the same as those in
Example 4, as detailed in Table~\ref{tab:setting}.
\begin{figure}[htb]
    \centering
    \includegraphics[width=0.82\linewidth]{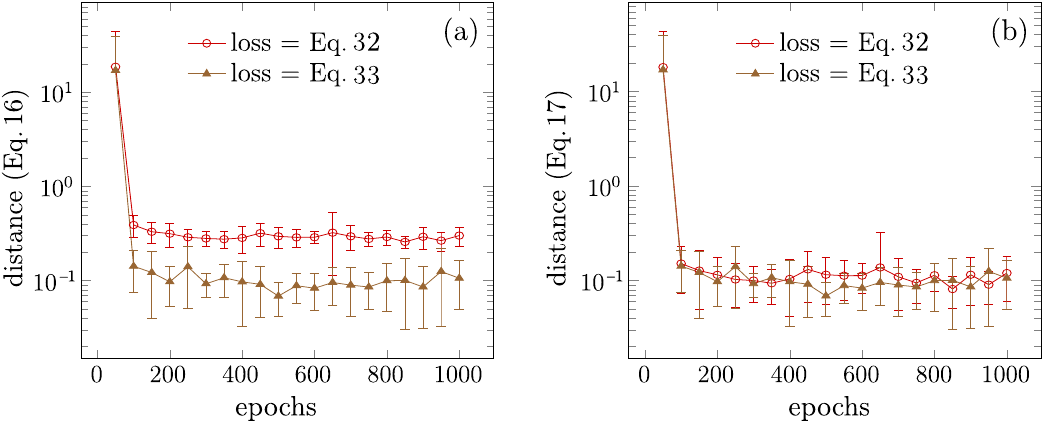}
    \caption{(a) The change in Eq.~\eqref{time_discretize} and
      Eq.~\eqref{approximation} when minimizing
      Eq.~\eqref{time_discretize} over training epochs. (b) The change
      in Eq.~\eqref{time_discretize} and Eq.~\eqref{approximation}
      when minimizing Eq.~\eqref{approximation} over training epochs.}
    \label{fig:comparison}
\end{figure}
First, we minimize
Eq.~\eqref{time_discretize} and record Eq.~\eqref{time_discretize} and
Eq.~\eqref{approximation} over training epochs. Next, we minimize
Eq.~\eqref{approximation} and record Eq.~\eqref{time_discretize} and
Eq.~\eqref{approximation} over training epochs. The results are shown
in Fig.~\ref{fig:comparison}.

From Fig.~\ref{fig:comparison} (a), we can see that when minimizing
Eq.~\eqref{time_discretize}, Eq.~\eqref{time_discretize} is almost
$10^{0.5}$ times larger than Eq.~\eqref{approximation}. However, when
minimizing Eq.~\eqref{approximation}, the values of
Eq.~\eqref{time_discretize} and Eq.~\eqref{approximation} are close to
each other (Fig.~\ref{fig:comparison} (b)). In both cases,
Eq.~\eqref{approximation} converges to approximately $10^{-1}$.
Interestingly, minimizing Eq.~\eqref{approximation} leads to a smaller
value of Eq.~\eqref{time_discretize}. This again implies that
minimizing Eq.~\eqref{approximation} can be more effective than
minimizing Eq.~\eqref{time_discretize} in Example~\ref{example4}. More
analysis on Eq.~\eqref{approximation} is needed to understand its
theoretical properties and to compare the performances of minimizing
Eq.~\eqref{approximation} versus minimizing
Eq.~\eqref{time_discretize} from numerical aspects.

\end{document}